\documentclass[a4paper,10pt,reqno]{amsart}
\usepackage[utf8]{inputenc}
\usepackage{amsmath, amssymb, mathrsfs}
\usepackage{bbm}
\usepackage{parskip}
\usepackage{mathtools}
\usepackage{csquotes}
\usepackage{verbatim}
\usepackage{bm}
\usepackage{ragged2e}
\usepackage{graphicx}
\usepackage{quiver}
\usepackage{amsthm}
\usepackage[hidelinks]{hyperref}
\usepackage[all]{xy}
\usepackage{ stmaryrd }
\usepackage{tikz}
\usepackage{tikz-cd}
\usepackage{ stmaryrd }
\usepackage{cleveref}
\usetikzlibrary{decorations.markings,positioning,calc}
\usepackage[style=alphabetic,sorting=nty,maxcitenames=50,maxnames=50]{biblatex}
\usepackage{color}   
\usepackage{hyperref}
\usepackage[left=3cm,right=3cm,top=3cm, bottom=4cm]{geometry}
\usepackage{url}

\hypersetup{
    colorlinks=true, 
    linktoc=all,     
    linkcolor=blue,  
}
\tikzset{
  closed/.style = {decoration = {markings, mark = at position 0.5 with { \node[transform shape, xscale = .8, yscale=.4] {/}; } }, postaction = {decorate} },
  open/.style = {decoration = {markings, mark = at position 0.5 with { \node[transform shape, scale = .7] {$\circ$}; } }, postaction = {decorate} }
}

\newcounter{diagram}  

\crefname{diagram}{diagram}{diagrams}
\crefname{diagram}{Diagram}{Diagrams}
\creflabelformat{diagram}{(#1) #2 #3}

\theoremstyle{definition}
\newtheorem{definition}{Definition}[section]
\newtheorem{remark}[definition]{Remark}
\newtheorem{example}[definition]{Example}

\theoremstyle{plain}
\newtheorem{lemma}[definition]{Lemma}
\newtheorem{theorem}[definition]{Theorem}
\newtheorem{proposition}[definition]{Proposition}
\newtheorem{corollary}[definition]{Corollary}

\newcommand{\overbar}[1]{\mkern 1.5mu\overline{\mkern-1.5mu#1\mkern-1.5mu}\mkern 1.5mu}

\DeclarePairedDelimiterX{\inp}[1]{\langle}{\rangle}{#1}

\def\id{\operatorname{id}}
\def\Spec{\operatorname{Spec}}
\def\relSpec{\underline{\Spec}}
\def\rs{\operatorname{rs}}
\def\rsInduced{\operatorname{rs}_*}
\def\rsInducedOverFnull{\operatorname{rs}_{*,f_0}}
\def\rsInducedRig{\operatorname{rs}_{*}^{\operatorname{rig}}}
\def\rsInducedFnullRig{\operatorname{rs}_{*,f_0}^{\operatorname{rig}}}
\def\Uhat{\hat{U}}
\def\evone{\operatorname{ev}_1}
\def\evnull{\operatorname{ev}_0}
\def\gr{\operatorname{gr}}
\def\Zrig{\mathfrak{Z}^{\operatorname{rig}}}
\def\grRig{\operatorname{gr}^{\operatorname{rig}}}
\def\plim{p_{\operatorname{lim}}}

\def\MapBGm{\underline{\operatorname{Map}}(B\mathbb{G}_m,\mathcal{X})}
\def\MapTheta{\underline{\operatorname{Map}}(\Theta,\mathcal{X})}

\def\Grad{\operatorname{Grad}(\mathcal{X})}
\def\Filt{\operatorname{Filt}(\mathcal{X})}
\def\FiltZ{\operatorname{Filt}(\mathcal{X})_{\mathfrak{Z}}}

\def\FiltZRig{\operatorname{Filt}(\mathcal{X})_{\mathfrak{Z}}\!\!\!\fatslash \mathcal{I}_{\operatorname{gr}}}

\def\FiltZnonsplit{\operatorname{Filt}(\mathcal{X})_{\mathfrak{Z}}^{\circ}}

\def\FiltZnonsplitRig{\operatorname{Filt}(\mathcal{X})_{\mathfrak{Z}}^{\circ}\!\!\!\fatslash \mathcal{I}_{\operatorname{gr}}}

\def\FiltQZ{\operatorname{Filt}([X/\hat{U}])_{Z\times B\mathbb{G}_m}}
\def\FiltQZRig{\operatorname{Filt}([X/\hat{U}])_{Z\times B\mathbb{G}_m}\!\!\!\fatslash \mathcal{I}_{\operatorname{gr}}}

\def\AttractZ{X^+_Z}
\def\AttractZnonsplit{X^{+,\circ}_Z}

\def\GmCan{\mathcal{G}^{\operatorname{can}}}

\newcommand{\Z}{\mathbb Z}

\newcommand{\coker}{\text{\coker }}

\newtheoremstyle{mytheorem}
{}
{}
{\itshape}
{}
{\scshape}
{:}
{.5em}
{}%
\newtheorem{innercustomthm}{Theorem}
\newenvironment{customthm}[1]
  {\renewcommand\theinnercustomthm{#1}\innercustomthm}
  {\endinnercustomthm}

\newtheoremstyle{mydefi}
{}
{}
{}
{}
{\scshape}
{:}
{.5em}
{}%

\theoremstyle{mytheorem}

\theoremstyle{mydefi}

\newtheorem{ex/}[teo]{Example}

\newtheorem{rmk/}[teo]{Remark}

\makeatletter
\renewenvironment{proof}[1][\proofname]{\par
  \vspace{-\topsep+6pt}
  \pushQED{\qed}%
  \normalfont
  \topsep0pt \partopsep0pt 
  \trivlist
  \item[\hskip\labelsep
        \scshape
    #1\@addpunct{.}]\ignorespaces
}{%
  \popQED\endtrivlist\@endpefalse
  \addvspace{6pt plus 6pt} 
}
\makeatother

\begin{document}
\title{Moduli spaces for \(\Theta\)-strata and non-reductive quotients}
\author{Ludvig Modin}
\address{Leibniz Universität Hannover, Institute für Algebraische Geometrie, Welfengarten 1, 30167 Hannover}
\email{modin@math.uni-hannover.de}

\begin{abstract}
    We give a new proof of the \(\Uhat\)-theorem  of Bérczi, Doran, Hawes and Kirwan on the existence of geometric quotients for actions of graded unipotent groups in terms of stacks of filtrations and gradings introduced by Halpern-Leistner. Our proof works over any affine Noetherian base, in particular it simultaneously generalizes the previous results to arbitrary characteristic, actions in families and to general \(\Theta\)-strata.
\end{abstract}

\maketitle
\section{Introduction}

Since algebraic stacks were introduced in \cite{Deligne-Mumford}, there has been a fruitful interplay between the study of their geometry and the equivariant geometry of schemes with a group action. This paper contributes to this interplay in two ways. First, we expand on the analogy between Białynicki-Birula strata of a scheme with a multiplicative action and \(\Theta\)-strata of an algebraic stack, introduced in \cite[Example 5.3.8]{halpernleistner2022structure} where Halpern-Leistner shows how to get a \(\Theta\)-stratification of a quotient stack \([X/\mathbb{G}_m]\) from the Białynicki-Birula stratification of \(X\). Secondly, we use the results obtained from the expansion of the first analogy to generalize one of the foundational results in non-reductive geometric invariant theory, the \(\Uhat\)-theorem. See \cite[Theorem 0.2]{U-hat} for the original theorem and \Cref{Uhat intro} for our generalization. Our proof of the \(\Uhat\)-theorem works over any affine Noetherian scheme, whereas the original proof works over an algebraically closed field of characteristic \(0\).

If \(X\) is smooth projective variety with an action of \(\mathbb{G}_m\), then Białynicki-Birula, \cite[Theorem 4.1]{BB}, associates to each fixed point component \(Z\subset X\), a locally closed subscheme \(X^+_Z\subset X\), the attractor of \(Z\), defined by \(x\in X_Z^+\) if and only if \({\lim_{t\to 0}t.x\in Z}\). He furthermore shows that \(X^+_Z\to Z\), defined by \(x\mapsto \lim_{t\to 0}t.x\) is affine and \(\mathbb{G}_m\)-equivariant. Together with the relative \(\operatorname{Proj}\)-construction, this allows for a construction of a geometric \(\mathbb{G}_m\)-quotient of \(X^+_Z\setminus Z\), projective over \(Z\). 

One can interpret \(Z\) and \(X^+_Z\) as subspaces of the equivariant mapping spaces \(\operatorname{Map}^{\mathbb{G}_m}(\operatorname{pt},X)\) and \(\operatorname{Map}^{\mathbb{G}_m}(\mathbb{A}^1,X)\) respectively, which was first done in \cite{Hesselink-concentration}. Under this identification, the map \(X^+_Z\to Z\) is given by restricting an equivariant map \(\mathbb{A}^1\to X\) to \(0\). 

The first step to going from Białynicki-Birula to \(\Theta\)-strata is replacing the equivariant mapping spaces above with the mapping stacks, \(\Grad=\operatorname{Map}(B\mathbb{G}_m,\mathcal{X})\) and \(\Filt=\operatorname{Map}([\mathbb{A}^1/\mathbb{G}_m],\mathcal{X})\) for an algebraic stack \(\mathcal{X}\). These are called the stacks of gradings and filtrations of \(\mathcal{X}\), and they get their name as they parameterize graded and filtered vector bundles respectively if \(\mathcal{X}\) is a stack of vector bundles on a curve. Halpern-Leistner shows in \cite{halpernleistner2022structure} that \(\Grad\) and \(\Filt\) are algebraic under reasonable assumptions on \(\mathcal{X}\). 

Restriction along \(B\mathbb{G}_m\subset [\mathbb{A}^1/\mathbb{G}_m]\) defines a morphisms \(\gr:\Filt\to \Grad\), analogous to the limit morphism to a fixed point component from its attractor. Similarly, restrictions along the inclusion of the open point \(\operatorname{pt}\subset [\mathbb{A}^1/\mathbb{G}_m]\) and the natural map \([\mathbb{A}^1/\mathbb{G}_m]\to B\mathbb{G}_m\) defines \(\evone:\Filt\to \mathcal{X}\) and \({\sigma:\Grad \to \Filt}\), respectively analogous to the inclusion of an attractor in the whole space and the inclusion of a fixed point component in its attractor. 

The fiber \(\FiltZ\) of \(\gr:\Filt\to \Grad\) over a substack \(\mathfrak{Z}\subset \Grad\)  is a \(\Theta\)-stratum of \(\mathcal{X}\) if \(\evone:\FiltZ\to \mathcal{X}\) is a locally closed immersion.


Since the limit map of a Białynicki-Birula stratum is affine, it is natural to ask if \({\gr:\FiltZ\to \mathfrak{Z}}\) is as well. In general this can not hold, as \({\gr:\Filt\to\Grad}\) doesn't have to be representable. 

The analogue of removing the fixed points of an attractor is to remove the section \(\sigma(\mathfrak{Z})\). We call \({\FiltZnonsplit=\FiltZ\setminus \sigma(\mathfrak{Z})}\) the stack of non-split filtrations centered at \(\mathfrak{Z}\). Based on this analogy, it is natural to ask if \(\FiltZnonsplit\) admits a coarse moduli space with some projectivity properties. This too is not true in general to my knowledge, but the answer to this question and the previous one are closely related. We show:

\begin{customthm}{A}\label{Theorem A}
    Let \(\mathcal{X}\) be an algebraic stack locally of finite presentation, with separated and quasi-compact diagonal and affine stabilizers over a quasi-separated scheme \(S\). Let \(\mathfrak{Z}\subset \Grad\) be a connected stack of gradings of \(\mathcal{X}\), then
    \begin{enumerate}
        \item there is a central flat subgroup stack \(\GmCan\subset \mathcal{I}_{\mathfrak{Z}}\) of the inertia stack of \(\mathfrak{Z}\), either isomorphic to \(\mathbb{G}_{m,\mathfrak{Z}}\) or the trivial group.
    \end{enumerate}
    Let \(\Zrig:=\mathfrak{Z}\!\!\!\fatslash \GmCan\) be the rigidification of \(\mathfrak{Z}\) with respect to \(\GmCan\).
    
    If the relative inertia of \(\FiltZ\to \mathfrak{Z}\), \(\mathcal{I}_{\gr}\), is flat of finite presentation over \(\FiltZ\), then
    \begin{enumerate}
        \setcounter{enumi}{1}
        \item the induced morphism \(\FiltZRig\to \mathfrak{Z}\),  from the rigidification w.r.t. \(\mathcal{I}_{\gr}\), is affine of finite presentation and,
        \item the substack \(\FiltZnonsplitRig\subset \FiltZRig\) is  open and the induced morphism \[\FiltZnonsplitRig\to \Zrig\] admits a tame relative moduli space, smooth locally projective over \(\mathfrak{Z}^{rig}\).
    \end{enumerate}
    In particular, if \(\gr:\FiltZ\to \mathfrak{Z}\) is representable, it is affine of finite presentation and \(\FiltZnonsplit\) admits a locally projective tame relative moduli space over \(\Zrig\).
\end{customthm}


The first part of this theorem is \Cref{gmcan properties}, the second and third are  \Cref{Part 2 of theorem A} and \Cref{Part 3 of theorem A} respectively. 

\begin{remark}
    We note some cases where part 3 is nicer:
    \begin{itemize}
        \item if \(\mathfrak{Z}=\Zrig\), then \(\gr:\Filt\to \mathfrak{Z}\) is an isomorphism, and \({\FiltZnonsplitRig}\) is empty,
        \item if \(\mathfrak{Z}\cong \Zrig\times B\mathbb{G}_m\), then the relative moduli space of \({\FiltZnonsplitRig\to \Zrig}\) is projective and,
        \item if  \(\Zrig\) is a Deligne-Mumford stack, then the relative moduli space of \({\FiltZnonsplitRig\to \Zrig}\) is étale locally projective.
    \end{itemize}
    The reason we can say more in the above cases is that to show projectivity of the relative moduli space of 
    \(\FiltZnonsplitRig\to \Zrig\), we need to neutralize the \(\mathbb{G}_m\)-gerbe \(\mathfrak{Z}\to \Zrig\). This can be done étale locally if \(\Zrig\) is a DM stack whereas it can only be done smooth locally if \(\Zrig\) is an Artin stack in general, as far as I am aware.
\end{remark}

Using \Cref{Theorem A}, we deduce a version of the \(\Uhat\)-theorem, \cite[Theorem 0.2]{U-hat}, on the existence of projective geometric quotients of certain non-reductive group actions on projective schemes. 

A smooth affine group scheme \(\Uhat\) over an affine scheme \(S\) is a graded unipotent group if 
\[\Uhat\cong U\rtimes \mathbb{G}_m \]
where \(U\) is a smooth unipotent group and the adjoint action of \(\mathbb{G}_m\) on \(\operatorname{Lie}(U)\) has strictly positive weights. Let \(\Uhat\) be a graded unipotent group acting on a projective scheme \(X\), linearized by an ample line bundle. Let \(Z\subset X\) be a fixed point component for the \(\mathbb{G}_m\) action on \(X\) induced by restriction along \(\mathbb{G}_m\subset \Uhat\), and let \(X^+_Z\) be the associated Białynicki-Birula stratum.

\begin{customthm}{B}\label{Uhat intro}
Let \(X\) be a reduced projective scheme with irreducible geometric fibers over a Noetherian affine scheme \(S\), with an action of a graded unipotent group \(\hat{U}\), linearized by a very ample line bundle.
Let \(Z\subset X\) be the fixed point component corresponding to the open Białynicki-Birula stratum.
Suppose the action of \(U\) on \(X_{Z}^+\) has smooth stabilizers with constant dimension, then 
\begin{enumerate}
    \item the action of \(U\) on \(X_{Z}^+\) admits a geometric quotient affine over \(Z\)
    \item the action of \(\hat{U}\) on \(X^+_Z\setminus \Uhat.Z\) admits a geometric quotient projective over \(S\).
\end{enumerate}
If the action of \(U\) on \(X_Z^+\) is free, the above quotients are tame, i.e. pushforward of equivariant quasi-coherent sheaves along the quotient morphism is exact. 
\end{customthm}
A slightly stronger version of this is \Cref{Classical Uhat} in the text. We deduce it from \Cref{Theorem B in text}, which applies to more general actions of graded unipotent groups, and naturally gives similar results for all Białynicki-Birula strata, not just the open one. All the results on formations of quotients in the paper are corollaries to \Cref{Theorem A}, using the description of the stacks of filtrations and gradings of quotient stacks in \cite[Theorem 1.4.8]{halpernleistner2022structure} over a field and \Cref{Filt of quotients} over a general Noetherian base scheme. 

In \Cref{comparison with classical uhat} we discuss in more detail how \Cref{Uhat intro} relates to the original \(\Uhat\)-theorem, in particular how different properties of the quotients relate to properties of the linearization of the action.

The main applications of \cite[Theorem 0.2]{U-hat} have been in developing geometric invariant theory (NRGIT) for (some) non-reductive group actions. This has mainly been carried out in \cite{U-hat}, \cite{bérczi2020projective}, \cite{hamilton2024affine}, \cite{hoskins2021quotients}, \cite{StratifyingQuotients} and \cite{MomentMapNRGIT}. This in turn has seen a applications to the Kobayashi conjecture in \cite{HyperbolicNRGIT}, using NRGIT to construct an alternative compactification of the bundle of regular jets, invariant under reparameterization. We intend to extend this compactification to invariant jets of schemes over a Noetherian base in future work.
Another application of NRGIT has been to construct moduli spaces of unstable objects for moduli problems that admits a GIT presentation, for example in \cite{hoskins2021quotients}. We are working with George Cooper on applying \Cref{Theorem A} to such constructions, for moduli problems that don't have a natural GIT presentation, for example moduli of Bridgeland unstable objects in a triangulated category. 

The strategy of proof of \Cref{Theorem A} is to introduce an action of \(\mathbb{G}_m\) on the fibers of the morphism \(\gr:\FiltZ\to \mathfrak{Z}\subset \Grad\) and applying a generalization of the Białynicki-Birula decomposition to actions of \(\mathbb{G}_m\) on algebraic spaces. This generalization is due to Drinfeld \cite{drinfeld2015algebraic} over a field and Richarz \cite{richarz2018spaces} over a base scheme, using a generalization of Sumihiro's theorem due to Alper, Hall and Rydh \cite{alper2023etale}. In particular, the proof of \Cref{Theorem A} is another instance of the interplay between the geometry of algebraic stacks and equivariant geometry.

The outline of the paper is as follows. In section 2, we recall the generalization of the Białynicki-Birula decomposition to algebraic spaces. We also give an explanation of what this says about quotient stacks of the form \([X/\mathbb{G}_m]\) to elucidate how \Cref{Theorem A} relates to the classical Białynicki-Birula decomposition.
In section 3, we recall some basic properties of \(\Theta\)-strata and discuss what they look like for quotient stacks and why it is natural to consider graded unipotent groups in this context.
In section 4 we prove \Cref{Theorem A}. In section 5 we prove \Cref{Uhat intro} and the other corollaries of \Cref{Theorem A} for quotient stacks. In the first appendix we prove a weakened version of \cite[Theorem 1.4.8]{halpernleistner2022structure} valid for quotient stacks over any Noetherian base. We need this to show \Cref{Uhat intro} when the base is not a field. In the second appendix we prove some formal properties of rigidifications of algebraic stacks, needed to prove \Cref{Theorem A}.

\subsection{Acknowledgments}
This paper is based on the first half of my PhD thesis, carried out in the DFG-Research Training Group 2553 at Universität Duisburg-Essen. The finalization this paper was carried out at Leibniz Universität Hannover.
I want to thank my supervisor Jochen Heinloth for his support throughout my thesis and for encouraging me to pursue \Cref{Theorem A} after I proved a version of \Cref{Uhat intro} directly. I thank Jochen Heinloth, again, and Jarod Alper for their comments on my thesis that helped improve the presentation of this paper.
I thank the Isaac Newton Institute and Victoria Hoskins for giving me the opportunity to present an early version of this project during the program "New equivariant methods in algebraic and differential geometry" that led to many fruitful discussions. Finally, I thank Wushi Goldring and Stefan Schreieder for useful comments on the first draft of this paper.
\section{Algebraic spaces with a multiplicative action}\label{gm-actions-on-spaces}
In this section we recall the generalization of the Białynicki-Birula decomposition to actions of \(\mathbb{G}_m\) on algebraic spaces due to Drinfeld \cite{drinfeld2015algebraic}, Halpern-Leistner \cite{halpernleistner2022structure}, Richarz \cite{richarz2018spaces} and Alper, Hall and Rydh \cite{alper2023etale}. We also give a a stacky reinterpretation of the classical Białynicki-Birula strata and explain how it relates to \Cref{Theorem A}.

\subsection{Equivariant mapping spaces}

Let \(S\) be a quasi-separated scheme, \(\mathbb{G}_m=\mathbb{G}_{m,S}\) the multiplicative group over \(S\) and \(\mathbb{A}^1=\mathbb{A}^1_S\) the affine line over \(S\) with the scaling action of \(\mathbb{G}_m\). That is, for any \(S\)-scheme \(T\), the action 
\[\mathbb{G}_m(T)\times \mathbb{A}^1(T)\to \mathbb{A}^1(T)\]
is defined by
\[(t,x)\mapsto tx\]
for all \(t\in \mathbb{G}_m(T)\) and all \(x\in \mathbb{A}^1(T)\). 

Let \(X\) be an algebraic space, quasi-separated and locally of finite presentation over \(S\), with an action of \(\mathbb{G}_m\).

Following Drinfeld, let
\[X^+:=\underline{\text{Map}}^{\mathbb{G}_m}(\mathbb{A}^1,X)\]
and
\[X^0:=\underline{\text{Map}}^{\mathbb{G}_m}(S,X)\]
be the sheaves of equivariant maps from \(\mathbb{A}^1\) respectively \(S\) into \(X\). Concretely, for any \(S\)-scheme \(T\), endowed with the trivial \(\mathbb{G}_m\)-action, 
\[X^+(T)=\{f:\mathbb{A}^1\times_ST\to X\mid f \text{ is }\mathbb{G}_m\text{-equivariant}\}\]
and
\[X^0(T)=\{f:T\to X\mid f \text{ is }\mathbb{G}_m\text{-equivariant}\}\]
in other words, \(X^0=X^{\mathbb{G}_m}\) parameterize the fixed points of the \(\mathbb{G}_m\)-action on \(X\).

These fit into the diagram

\begin{equation}\label{Situation 1}
\begin{tikzcd}
	{X^{+}} & X \\
	{X^0} & X
	\arrow["{\evone}", from=1-1, to=1-2]
	\arrow["{\plim}"', from=1-1, to=2-1]
	\arrow["{\evnull}"', from=2-1, to=2-2]
\end{tikzcd}
\end{equation}
with morphisms defined in the following way on \(T\)-points.
Let \(f\in X^+(T)\), then \(\evone(f):=f(1)\) and \(\plim(f)=f|_{\{0\}\times T}\).
Let \(g\in X^0(T)\), then \(\evnull(g)=g\) for \(g\in X^0(T)\), i.e. \(g:T\to X\) is an equivariant morphism, in particular it is a \(T\)-point of \(X\).
There is also a natural section \(s_0:X^0\to X^+\)
of \(\plim:X^+\to X^0\) defined by precomposing \(T\to X\) with the structure morphism \(\mathbb{A}^1_T\to T\), i.e. it takes fixed point to the corresponding constant morphism \(\mathbb{A}^1\to X\). 

Furthermore, \(\mathbb{A}^1\) is a monoid scheme under multiplication with a natural action on \(X^+\). Letting \(m:\mathbb{A}^1\times \mathbb{A}^1\to \mathbb{A}^1\) denote the multiplication on \(\mathbb{A}^1\), the action 
\[m_*:\mathbb{A}^1\times X^+\to X^+\]
is defined for \((t,f)\in \left(\mathbb{A}^1\times X^+\right)(T)\) by letting 
\[m_*(t,f)=f\circ m\circ (t\times \id_{\mathbb{A}^1})\]
i.e. we precompose \(f\) by multiplication by \(t\) on \(\mathbb{A}^1\).

\begin{definition}\label{contracting def}
    A \(\mathbb{G}_m\)-action on an algebraic space \(X\) is \emph{contracting} if \(\evone:X^+\to X\) is an isomorphism.
    If \(X\) is an algebraic space with a contracting action, then we identify \(X\) with \(X^+\) and denote by
    \[\plim:X\to X^0.\]
    the composition of this identification with \(\plim:X^+\to X^0\).
\end{definition}

\begin{example}\label{affine multip-act}
    Let \(X=\Spec \mathcal{A}\) be an affine \(S\)-scheme with an action given by a grading \({\mathcal{A}\cong \bigoplus_{i\in \Z}\mathcal{A}_i}\).
    Let \(\mathcal{I}_+\) be the ideal generated by \(\bigoplus_{i<0}\mathcal{A}_i\) and \(\mathcal{I}_{-}\) be the ideal generated by \(\bigoplus_{i>0}\mathcal{A}_i\), then \(X^+=\Spec \mathcal{A}/\mathcal{I}_+\) and \(X^0=\Spec \mathcal{A}/\left(\mathcal{I}_++\mathcal{I}_-\right)\), see \cite[Lemma 1.9]{richarz2018spaces}. In particular, the action of \(\mathbb{G}_m\) on \(X\) is contracting if and only if \(\mathcal{A}_i=0\) for all \(i<0\).
\end{example}

The key properties of multiplicative actions on algebraic spaces is summarized in the following theorem.

\begin{theorem}\label{multip-act}[\cite{drinfeld2015algebraic},\cite{richarz2018spaces},\cite{alper2023etale},\cite{halpernleistner2022structure}]
    Let \(X\) be an algebraic space quasi-separated and locally of finite presentation over \(S\) with an action of \(\mathbb{G}_m\).
    The following holds:
    \begin{enumerate}
        \item the morphism \(\evnull:X^0\to X\) is a closed immersion, in particular \(X^0\) is and algebraic space,
        \item the morphism \(\plim:X^+\to X^0\) is affine of finite presentation, in particular \(X^+\) is and algebraic space,
        \item the morphism \(\evone:X^+\to X\) is unramified,
         \item the fiber of \(\evone:X^+\to X\) over any geometric point of \(X^0\subset X\) has a single point,
        \item the \(\mathbb{G}_m\)-action on \(X\) is contracting if and only if 
        it extends to a monoid action of \((\mathbb{A}^1,\cdot)\) on \(X\), this extension is unique if it exists,
        \item if \(X\) is separated then \(\evone\) is a monomorphism,
        \item if \(X\) is proper then \(\evone\) induces a bijection on geometric points.
   
    \end{enumerate}
\end{theorem}

\begin{proof}
    Properties 1 and 2 is \cite[Theorem A]{richarz2018spaces}, which gives the result for étale locally linearizable actions, combined with \cite[Corollary 10.2]{alper2023etale}, which states that \(\mathbb{G}_m\)-actions on algebraic spaces \(X\), quasi-separated and locally of finite presentation over a quasi-separated base are étale locally linearizable.
    
    Property 6 is shown in \cite[Corollary 1.4.3]{halpernleistner2022structure}.
    
    The proofs of the remaining properties are given in the case \(S=\Spec k\) in \cite{drinfeld2015algebraic}. The proofs there carry over more or less identically to the case when \(S\) is any quasi-separated scheme once properties 1 and 2 have been established. For the interested reader, the details of this can be found in \cite[Theorem I 1.0.1]{ValCritUnstabStrat}.  
\end{proof}

\subsection{A stacky view of Białynicki-Birula strata}\label{Stacky BB}

    In this section we explain how much of the structure of Białynicki-Birula strata of a projective variety with an action of \(\mathbb{G}_m\) can be explained in terms of \(\Theta\)-strata and \Cref{Theorem A}. 
    
    Let \(X\) be a smooth projective variety over a field \(k\) with an action of \(\mathbb{G}_m\). 
    Let \(F_i\subset X^0\) be a connected component of the fixed point subscheme and \(X^+_i:=\plim^{-1}(F_i)\) its associated attractor. 
    By \cite[Theorem 5.27 (3c)]{Luna-slice-stacks}, the restriction of the evaluation map
    \[\evone:X^+\to X\]
    to \(X_i^+\) is a locally closed immersion . 

    The map \(\plim:X^+_i\to F_i\) is a good \(\mathbb{G}_m\) quotient, i.e \(F_i\) parametrize closed orbits of the action of \(\mathbb{G}_m\) on \(X^+_i\), in particular, if the action is not trivial, then this is not a geometric quotient. One way to get a geometric quotient is to remove \(F_i\) from \(X_i^+\). 
    Indeed, it follows from \(X_i^+\to F_i\) being affine and \(\mathbb{G}_m\)-equivariant that \(X_i^+\cong \underline{\Spec }_{\mathcal{O}_{F_i}}\mathcal{A}\) for \(\mathcal{A}=\bigoplus_{i\geq 0}\mathcal{A}_i\) a graded \(\mathcal{O}_{F_i}\)-algebra, with positive grading as the \(\mathbb{G}_m\)-action on \(X_i^+\) is contracting. Furthermore, the inclusion \(F_i\subset X_i^{+}\) is defined by the ideal \(\mathcal{I}=\bigoplus_{i> 0}\mathcal{A}_i\), the irrelevant ideal of \(\mathcal{A}\) in the construction of \(\underline{\operatorname{Proj}}_{\mathcal{O}_{F_i}}\mathcal{A}\). It follows from \cite[Proposition 1.6.1 (iii)]{WeightedBlowUps} that \(X_i^+\setminus F_i\to \underline{\operatorname{Proj}}_{\mathcal{O}_{F_i}}\mathcal{A} \) is a geometric quotient. All of this can be summarized in the following diagram of quotient stacks and their moduli spaces:
    
    \begin{equation}\label{Situation 2}
    \begin{tikzcd}
	{[(X_i\setminus F_i)/\mathbb{G}_m]} && {[X_{F_i}/\mathbb{G}_m]} && X \\
	{(X_i\setminus F_i)/\mathbb{G}_m} && {F_i\times B\mathbb{G}_m} \\
	&& {F_i}\tikzset{
  closed/.style = {decoration = {markings, mark = at position 0.5 with { \node[transform shape, xscale = .8, yscale=.4] {/}; } }, postaction = {decorate} },
  open/.style = {decoration = {markings, mark = at position 0.5 with { \node[transform shape, scale = .7] {$\circ$}; } }, postaction = {decorate} }
}
	\arrow[hook,open, from=1-1, to=1-3]
	\arrow["\text{cms}", from=1-1, to=2-1]
	\arrow["{\text{loc. closed}}", hook, from=1-3, to=1-5]
	\arrow["\text{affine}", from=1-3, to=2-3]
	\arrow["\text{projective}"', from=2-1, to=3-3]
	\arrow["{\mathbb{G}_m\text{-gerbe}}", from=2-3, to=3-3]
\end{tikzcd}
\end{equation}
where "cms" stands for coarse moduli space.

To relate this to \(\Theta\)-strata, we note that by \cite[Theorem 1.4.7]{halpernleistner2022structure}, when \({\mathcal{X}=[X/\mathbb{G}_m]}\) for the above action, the following diagram 
\[\begin{tikzcd}
	{[X_i^+/\mathbb{G}_m]} & \Filt \\
	{F_i\times B\mathbb{G}_m} & \Grad
	\arrow[from=1-1, to=1-2]
	\arrow["\plim"', from=1-1, to=2-1]
	\arrow["\gr", from=1-2, to=2-2]
	\arrow[from=2-1, to=2-2]
\end{tikzcd}\]
is cartesian with horizontal arrows, defined by descending equivariant maps \(\mathbb{A}^1\to X\) and \({\Spec k\to X}\) to maps of quotient stacks \([\mathbb{A}^1/\mathbb{G}_m]\to [X/\mathbb{G}_m]\) and \(B\mathbb{G}_m\to [X/\mathbb{G}_m]\). Moreover, the horizontal arrows are clopen immersions, \(\sigma:\Grad\to \Filt\) restricts to the closed immersion \(F_i\times B\mathbb{G}_m\to [X^+_i/\mathbb{G}_m]\) and letting \(\mathfrak{Z}=F_i\times B\mathbb{G}_m\subset \Grad\) we have \(\Zrig=F_i\) by \Cref{Gm-can/RigCent example}. 

If the action of \(\mathbb{G}_m\) on \(X\) is linearized by an ample line bundle, then \cite[Example 5.3.8]{halpernleistner2022structure} implies that \(\evone:[X^+_i/\mathbb{G}_m]\to [X/\mathbb{G}_m]\) is a locally closed immersion, i.e that \([X^+_i/\mathbb{G}_m]\) is a \(\Theta\)-stratum. 
The rest of the properties of (\ref{Situation 2}) follows from \Cref{Theorem A} once we note that \[{\gr|_{[X^+_i/\mathbb{G}_m]}:[X_i^+/\mathbb{G}_m]\to F_i\times B\mathbb{G}_m}\] is representable, except we only get that \((X_i\setminus F_i)/\mathbb{G}_m\to F_i\) is locally projective. To get projectivity we may apply \Cref{Theorem A when Zrig space}.

\section{$\Theta$-strata and graded unipotent groups}
Here we recall the basic definitions and properties of \(\Theta\)-strata from \cite{halpernleistner2022structure} needed to set up and prove \Cref{Theorem A}. We then recall what the stacks of filtrations and gradings look like in the case of quotient stacks. Finally, we discuss graded unipotent groups, and why they are natural to study from the point of view of \(\Theta\)-strata.

\subsection{Filtrations and gradings of algebraic stacks}

Let \(S\) be a quasi-separated scheme and let \(\mathbb{G}_m=\mathbb{G}_{m,S}\) act on \(\mathbb{A}^1=\mathbb{A}^1_S\) by scaling as in the previous section, that is the action is given by the morphism 

\[\mathbb{G}_m\times\mathbb{A}^1\to \mathbb{A}^1\]
defined on \(T\)-points for any \(S\)-scheme \(T\) by
\[(t,x)\mapsto tx.\]
Let 
\[\Theta:= [\mathbb{A}^1/\mathbb{G}_m]\]
be the quotient stack defined by this action.

By \cite[Proposition 10.3.7]{OlssonSpaSta}, maps \(T\to \Theta\) corresponds to pairs \((\mathcal{L},s:\mathcal{L}\to \mathcal{O}_T)\) where \(\mathcal{L}\) is a line bundle on \(T\) and \(s\) is a morphism of \(\mathcal{O}_T\)-modules. This can also be seen directly as follows. By the definition of quotient stacks, maps \(T\to \Theta\) corresponds to equivariant maps \(f:\mathcal{P}\to \mathbb{A}^1_T\), where \(\mathcal{P}\) is a  \(\mathbb{G}_m\)-torsor over \(T\). The contracted product \(\mathcal{P}\times^{\mathbb{G}_m}\mathbb{A}^1_T\) yields a line bundle over \(T\), and the equivariant map \(f:\mathcal{P}\to \mathbb{A}^1_T\) extends to a unique equivariant map \(\overbar{f}:\mathcal{P}\times^{\mathbb{G}_m}\mathbb{A}^1_T\to \mathbb{A}^1_T\). Let \(\mathcal{L}\) denote the sheaf of sections of \(\mathcal{P}\times^{\mathbb{G}_m}\mathbb{A}^1_T\). The map \(\overbar{f}:\mathcal{P}\times^{\mathbb{G}_m}\mathbb{A}^1_T\to \mathbb{A}^1_T\) defines a morphism of \(\mathcal{O}_T\)-modules, \(\mathcal{L}\to \mathcal{O}_T\) by \(s\mapsto \overbar{f}\circ s\), for \(s:T\to \mathcal{P}\times^{\mathbb{G}_m}\mathbb{A}^1_T\). Conversely, given a line bundle \(\mathcal{L}\) with a  map \(\mathcal{L}\to \mathcal{O}_T\), we get a \(\mathbb{G}_m\)-torsor with an equivariant map to \(\mathbb{A}^1_T\) by restricting to the complement of the zero section.
Similarly, \(B\mathbb{G}_m(T)\) is equivalent to the groupoid of line bundles over \(T\).

Define 
\[0:B\mathbb{G}_m\to \Theta\] 
by mapping mapping a line bundle \(\mathcal{L}\) over an \(S\)-scheme \(T\) to \((\mathcal{L},0)\),
\[q:\Theta\to B\mathbb{G}_m\] 
by mapping a pair \((\mathcal{L},s)\) over an \(S\)-scheme \(T\) to \(\mathcal{L}\), and 
\[1:S\to \Theta\]
by the pair \((\mathcal{O}_S,1)\).

Given an algebraic stack \(\mathcal{X}\) over \(S\), let

\begin{equation}\label{def of filt}
    \Filt:=\MapTheta
\end{equation}
and
\begin{equation}\label{def of grad}
        \Grad:=\MapBGm  
\end{equation}
be the stacks of maps \(\Theta\to \mathcal{X}\) and \(B\mathbb{G}_m\to \mathcal{X}\) respectively.
Over any \(S\)-scheme \(T\) these stacks have the groupoids of morphisms \(\Theta_T\to \mathcal{X}\) and \(B\mathbb{G}_{m,T}\to \mathcal{X}\) as fibers respectively.


These mapping stacks fit in a diagram similar to \ref{Situation 1},
\[\begin{tikzcd}\label{Situation 3}
	{\Filt} & {\mathcal{X}} \\
	{\Grad} & {\mathcal{X}}
	\arrow["{\gr}"', from=1-1, to=2-1]
	\arrow["{\evone}", from=1-1, to=1-2]
	\arrow["u"', from=2-1, to=2-2]
\end{tikzcd}\]
where for any \(S\)-scheme \(T\), these maps are defined as follows on \(T\)-points \(f\in \Filt\) and \({g\in \Grad}\):
\begin{align}\label{def of evone}
\begin{split}
  \evone:\Filt\to \mathcal{X}\\
  f\mapsto 1\circ f 
  \end{split}
\end{align}
\begin{align}
\begin{split}
\gr:\Filt\to \Grad\\
 f\mapsto 0\circ f
 \end{split}
 \end{align}
and 
\begin{align}
    \begin{split}
        u:\Grad\to \mathcal{X}\\
        g\mapsto p\circ g
    \end{split}
\end{align}
where \(p:T\to B\mathbb{G}_{m,T}\) is the atlas of defined by the trivial action on \(T\), i.e it corresponds to the trivial line bundle \(\mathcal{O}_T\) on \(T\).

There is also a section of \(\gr:\Filt\to \Grad\), 
\begin{align}
    \begin{split}
        \sigma:\Grad\to \Filt\\
        g\mapsto g\circ q
    \end{split}
\end{align}
just like there is an inclusion of \(\mathbb{G}_m\)-fixed points into the attractor of a scheme with a \(\mathbb{G}_m\)-action.

We summarize some key properties of these stacks shown in \cite[Propositions 1.1.2 and 1.3.8]{halpernleistner2022structure}:

\begin{theorem}\cite{halpernleistner2022structure}\label{basic-prop-of-filt-grad}
    Let \(\mathcal{X}\) be an algebraic stack locally of finite presentation, with separated and quasi-compact diagonal and affine stabilizers over a quasi-separated scheme \(S\). Then the following holds
    \begin{enumerate}
        \item the stacks \(\Filt\) and \(\Grad\) are algebraic stacks with separated  and quasi-compact diagonal with affine stabilizers, locally of finite presentation over \(S\),
        \item the morphism \(\gr:\Filt\to \Grad\) is of finite presentation and induces a bijection on connected components.
    \end{enumerate}
\end{theorem}
\begin{proof}
    The first statements is shown in \cite[Proposition 1.1.2]{halpernleistner2022structure}.
    In \cite[Proposition 1.3.8]{halpernleistner2022structure}, it is shown that \(\gr:\Filt\to \Grad\) is quasi-compact and induces a bijection of connected components with inverse induced by \(\sigma\). By part 1, both \(\Filt\) and \(\Grad\) are quasi-separated and locally of finite presentation over \(S\), hence \(\gr:\Filt\to \Grad\) is a quasi-separated morphism locally of finite presentation. Because it is also quasi-compact, it is of finite presentation.  
\end{proof}

We will be interested in the fibers of \(\gr:\Filt\to \Grad\), for this we introduce the following notation.

\begin{definition}\label{filt-centered-at-Z}
    Given a morphism of \(S\)-stacks \(\mathfrak{Z}\to \Grad)\), we call the pullback 
    \[\begin{tikzcd}
	{\FiltZ} & {\Filt} \\
	{\mathfrak{Z}} & {\Grad}
	\arrow[from=1-1, to=1-2]
	\arrow[from=1-2, to=2-2]
	\arrow[from=1-1, to=2-1]
	\arrow[from=2-1, to=2-2]
\end{tikzcd}\]
the \textit{stack of filtrations centered at \(\mathfrak{Z}\)}. 

Let \(\sigma(Z)\subset \FiltZ\) be the essential image of \(\sigma:\mathfrak{Z}\to \FiltZ\) viewed as a morphism of fibered categories and let
\[\FiltZnonsplit:=\FiltZ\setminus \sigma(\mathfrak{Z})\]
we call this the \emph{stack of non-split filtrations centered at \(\mathfrak{Z}\)}.
\end{definition}
\begin{remark}
    In the setting of \Cref{Theorem A}, \(\sigma(Z)\) is a closed substack, as we show in \Cref{Part 2 of theorem A}, but in general it is not clear that the essential image of \(\sigma:\mathfrak{Z}\to \FiltZ\) agrees with the scheme theoretic image.
\end{remark}

A special case of these stacks of filtrations are \(\Theta\)-strata of \(\mathcal{X}\)
\begin{definition}\label{theta-stratum}\cite[Definition 2.5.3]{núñez2023refinedhardernarasimhanfiltrationsmoduli}
    The stack of filtrations centered at \(\mathfrak{Z}\subset \Grad\), \(\FiltZ\), is a \textit{\(\Theta\)-stratum} of \(\mathcal{X}\) if
    the morphism 
    \[\evone|_{\FiltZ}:\FiltZ\to \mathcal{X}\] 
    is a locally closed immersion. The stack \(\mathfrak{Z}\subset \mathcal{X}\) is called the \textit{center} of the \(\Theta\)-stratum \(\FiltZ\).
\end{definition}

This version of the definition of a \(\Theta\)-stratum, due to Ibáñez Núñez, is a bit more flexible than the original definition \cite[Definition 2.1.1]{halpernleistner2022structure}, which requires \(\FiltZ\) to be a union of components and \(\evone|_{\FiltZ}:\FiltZ \to \mathcal{X}\) to be a closed immersion.

\subsection{Filtrations and gradings of quotient stacks}\label{filt of quotients in text}
Here we state the known results on when filtrations an gradings of quotient stacks admit an explicit quotient stack description.

Let \(S\) be a scheme and let \(G\) be an affine group scheme over \(S\) acting on an algebraic space \(X\) quasi-separated and locally of finite presentation over \(S\), and let \(\lambda:\mathbb{G}_m\to G\) be a cocharacter. The action of \(G\) on \(X\) restricts via \(\lambda\) to an action of \(\mathbb{G}_m\) on \(X\), let \(X^{+,\lambda}\) and \(X^{0,\lambda}\) be the attractor and fixed points of the \(\mathbb{G}_m\)-action on \(X\).
The composition of \(\lambda\) with the conjugation action on \(G\) induces an action of \(\mathbb{G}_m\) on \(G\), let \(P_{\lambda}\) and \(L_{\lambda}\) be the attractor and fixed points of the action of \(\mathbb{G}_m\) on \(G\) and let \(U_{\lambda}\) be the kernel of \(\plim:P_{\lambda}\to L_{\lambda}\). Note that \(L_{\lambda}\) is the centralizer in \(G\) of \(\lambda\).

\begin{remark}
    If \(G\) is a reductive algebraic group, then \(P_{\lambda}\) is the parabolic subgroup associated to \(\lambda\), \(L_{\lambda}\) is the associated Levi subgroup and \(U_{\lambda}\) the unipotent radical of \(P_{\lambda}\).
\end{remark}

Since the functors \((-)^+\) and \((-)^0\) commute with products, there are natural group actions of \(P_{\lambda}\) on \(X^{+,\lambda}\) and of \(L_{\lambda}\) on \(X^0\). 

Let 
\[F_{\lambda}:{[X^{+,\lambda}/P_{\lambda}]} \to  {\operatorname{Filt}([X/G])}\]
be the morphism defined by mapping \(f\in X^{+,\lambda}(T)\) to the induced morphism \(\Theta_T\to [X/G]\), 
\[G_{\lambda}:[X^{0,\lambda}/L_{\lambda}] \to \operatorname{Grad}([X/G])\]
the morphism defined by mapping \(g\in X^{0,\lambda}(T)\) to the induced morphism \(B\mathbb{G}_{m,T}\to [X/G]\) and
\[P_{\lim}:[X^{+,\lambda}/P_{\lambda}]\to [X^{0,\lambda}/L_{\lambda}]\] defined by \(\plim:X^{+,\lambda}\to X^{0,\lambda}\) and the quotient morphism \(P_{\lambda}\to P_{\lambda}/U_{\lambda}=L_{\lambda}\). 

\begin{theorem}[Theorem 1.4.7 in \cite{halpernleistner2022structure}]
Let \(X\) be an algebraic space with an action of \(\operatorname{GL}_n\) and let \(\mathcal{X}=[X/\operatorname{GL}_n]\). 
Let \(\operatorname{Hom}(\mathbb{G}_m,T)\) denote the cocharacter group of the maximal torus of diagonal matrices with the natural action of \(S_n\). 
The diagram
\[\begin{tikzcd}
	{\coprod_{\lambda\in \text{Hom}(\mathbb{G}_m,T)/S_n}[X^{+,\lambda}/P_{\lambda}]} & \Filt \\
	{\coprod_{\lambda\in \text{Hom}(\mathbb{G}_m,T)/S_n}[X^{0,\lambda}/L_{\lambda}]} & \Grad
	\arrow["\cong",from=1-1, to=1-2]
	\arrow["P_{\lim}"', from=1-1, to=2-1]
	\arrow["\gr", from=1-2, to=2-2]
	\arrow["\cong", from=2-1, to=2-2]
\end{tikzcd}\]
is cartesian, with horizontal arrows given by \(F_{\lambda}:{[X^{+,\lambda}/P_{\lambda}]} \to  \Filt\) and \(G_{\lambda}:[X^{0,\lambda}/L_{\lambda}] \to \operatorname{Grad}([X/G])\) for each fixed \(\lambda\).
\end{theorem}
In \cite[Theorem 1.4.8]{halpernleistner2022structure}, a similar result is proven for actions of smooth affine group schemes over a field.
For our applications of \Cref{Theorem A} to non-reductive quotients, we will need a version of the above for more general group schemes over any Noetherian base, but we don't need to fully describe \(\Filt\) and \(\Grad\), only some substacks. We prove the following in the first appendix.

\begin{theorem}[\Cref{Filt of quotients}]
    Let \(G\) be a smooth connected affine group scheme over a Noetherian scheme \(S\) with an action on a quasi-separated algebraic space \(X\) locally of finite type over \(S\) and let \(\lambda:\mathbb{G}_{m}\to G\)  be a cocharacter. The diagram 
    \[\begin{tikzcd}
	{[X^{+,\lambda}/P_{\lambda}]} & {\operatorname{Filt}([X/G])} \\
	{[X^{0,\lambda}/L_{\lambda}]} & {\operatorname{Grad}([X/G])}
	\arrow["F_{\lambda}", from=1-1, to=1-2]
	\arrow["P_{\lim}"',from=1-1, to=2-1]
	\arrow["G_{\lambda}",from=2-1, to=2-2]
	\arrow["\gr",from=1-2, to=2-2]
\end{tikzcd}\]
commutes, and the horizontal arrows are open immersions.
\end{theorem}

We will also need the following lemma, which is a consequence of the definitions of 

\(F_{\lambda}:[X^{+,\lambda}/P_{\lambda}]\to\operatorname{Filt}([X/G])\) and \(G_{\lambda}:[X^{0,\lambda}/L_{\lambda}]\to \operatorname{Grad}([X/G])\).

\begin{lemma}\label{Non-split filt in quotients}
Let 
\[i_{\lambda}:[Z/L_{\lambda}]\to [X^{+,\lambda}/P_{\lambda}]\]
be the morphism induced by the inclusions \(Z\subset X^{0,\lambda}\subset X^{+,\lambda}\) and \(L_{\lambda}\subset P_{\lambda}\). 
Then \[F_{\lambda}\circ i_{\lambda}=\sigma \circ G_{\lambda}.\]

In particular, \(\FiltZnonsplit=[(\plim^{-1}(Z)\setminus P_{\lambda}.Z)/P_\lambda]\).
\end{lemma}
\begin{proof}
    Consider the commutative diagram 
    \[\begin{tikzcd}
	{X^{+,\lambda}} & {\operatorname{Filt}([X/G])} \\
	{X^{0,\lambda}} & {\operatorname{Grad}([X/G])}
	\arrow[from=1-1, to=1-2]
	\arrow["\plim"', from=1-1, to=2-1]
	\arrow["\gr", from=1-2, to=2-2]
	\arrow[from=2-1, to=2-2]
\end{tikzcd}\]
with horizontal arrows defined by mapping a \(\lambda(\mathbb{G}_m)\)-equivariant map \(\mathbb{A}^1\to X\) to the corresponding \(\Theta\to [X/G]\) defined by descent and similarly for the lower horizontal arrow. In other words, the horizontal arrows factor through \(F_{\lambda}\) and \(G_{\lambda}\). Now, the inclusion \(X^{0,\lambda}\subset X^{+,\lambda}\) corresponds to pre-composing \(\lambda(\mathbb{G}_m)\)-equivariant maps \(T\to X \) with the structure map \(\mathbb{A}^1_T\to T\). This descends to the pre-composition of the induced map \(B\mathbb{G}_m\to [X/G]\) with \(\Theta_T\to B\mathbb{G}_{m,T}\). In other words, \(F_{\lambda}\circ i_{\lambda}=\sigma\circ G_{\lambda}\), proving the lemma.
\end{proof}

\subsection{Graded unipotent groups}
With \cite[Theorems 1.4.7 and 1.4.8]{halpernleistner2022structure} and \Cref{Filt of quotients} in mind, it is natural to consider groups admitting cocharacters \(\lambda\) such that \(P_{\lambda}\) and \(L_{\lambda}\) are particularly well behaved. It is for actions of such groups that \Cref{Theorem A} can be applied directly to prove existence of quotients as we do in \Cref{Uhat intro}.

\begin{definition}\label{graded-unip-def}
    A smooth affine group scheme \(\Uhat\) of finite presentation over \(S\) is a \emph{graded unipotent group} if there is a cocharacter \(\lambda:\mathbb{G}_m\to \Uhat\), such that:    
    \begin{enumerate}
        \item \(P_{\lambda}=\Uhat\),
        \item \(L_{\lambda}=\mathbb{G}_m\) and
        \item \(U_{\lambda}=\ker(\plim: P_{\lambda}\to L_{\lambda})\) is smooth with unipotent geometric fibers over the base \(S\).
    \end{enumerate}
In this case we say that \(U_{\lambda}\) is \emph{graded} by the \(\mathbb{G}_m\)-action through conjugation via \(\lambda:\mathbb{G}_m\to \Uhat\).
\end{definition}

\begin{remark}
Part (3) of the definition is automatic, we only include it for clarity. Indeed, note that each geometric fiber of \(P_{\lambda}\) is a smooth affine variety over an algebraically closed field with a contracting \(\mathbb{G}_m\)-action, and the fiber of \(U_{\lambda}\) is the fiber over a point of the limit map \(P_{\lambda}\to L_{\lambda}\). Hence the fibers are isomorphic as schemes to \(\mathbb{A}^n\) for some \(n\geq 0\) by \cite[Theorem 4.1]{BB}. By Lazard's theorem, see e.g \cite[§4, 4.1(Thm. of Lazard)]{DemazureGabriel}, any algebraic group over an algebraically closed field, isomorphic to \(\mathbb{A}^n\) as a scheme, is unipotent.
\end{remark}

\begin{remark}\label{comparisson-of-graded-unipotent-defs}
    In \cite[Definition 2.1]{U-hat}, a complex unipotent group \(U\) is said to be graded by a \(\mathbb{G}_m\)-action through homomorphisms if the induced action on \(\mathfrak{u}=\text{Lie}(U)\) is of the form
    \[\mathfrak{u}=\bigoplus_{i>0} \mathfrak{u}_i\]
    where \(\mathbb{G}_m\) acts with weight \(i\) on the vector space \(\mathfrak{u}_i\).
    We claim that this definition coincides with \Cref{graded-unip-def} for complex groups.
    Note that if \(U\) is graded by a \(\mathbb{G}_m\)-action in the sense of \cite{U-hat}, the adjoint action on \(\mathfrak{u}\) is contracting. In characteristic \(0\), the exponential map \(\exp: \mathfrak{u}\to U\) is an isomorphism of varieties if \(U\) is unipotent, and by the Baker-Campell-Hausdorff formula, it is also \(\mathbb{G}_m\)-equivariant, since the action respects the Lie-bracket.
    Hence the \(\mathbb{G}_m\)-action on \(U\) also extends to an \(\mathbb{A}^1\)-action via the isomorphism of varieties \(\exp:\mathfrak{u}\to U\). 
    Conversely, if \(U\) is graded by a \(\mathbb{G}_m\)-action in the sense of \Cref{graded-unip-def}. The inverse of \(\exp:\mathfrak{u}\to U\), \(\log:U\to \mathfrak{u}\) is equivariant by the Baker-Campell-Hausdorff formula, hence the linear \(\mathbb{G}_m\)-action on \(\mathfrak{u}\) is contracting with only \(0\) as a fixed point, so after diagonalizing we get that
    \[\mathfrak{u}=\bigoplus_{i>0} \mathfrak{u}_i\]
    where \(\mathbb{G}_m\) acts with weight \(i\) on \(\mathfrak{u}_i\).
\end{remark}
\begin{lemma}\label{semi-dir decomp}
    Any graded unipotent group \(\Uhat\) over a scheme \(S\) admits a semi-direct product decomposition 
    \[\Uhat=U\rtimes \mathbb{G}_m\]
    where \(U\) is an affine group scheme over \(S\) with unipotent geometric fibers.
\end{lemma}
\begin{proof}
As \(\Uhat\) is graded unipotent, there is a cocharacter \(\lambda:\mathbb{G}_m\to \Uhat\) such that the conjugation action is contracting and the fixed point subscheme is isomorphic to \(\mathbb{G}_m\), which is a subgroup since the \(\mathbb{G}_m\)-action is by homomorphisms. Moreover, the inclusion \(\mathbb{G}_m\subset \Uhat\) is a section of the limit morphism \(\Uhat\to \mathbb{G}_m\). Finally, again since the action is via homomorphisms, \(\Uhat\to \mathbb{G}_m\) is a homomorphism, so the section \(\mathbb{G}_m\subset \Uhat\) defines a semi-direct product decomposition.
\end{proof}

A basic example of a graded unipotent group to keep in mind is the following. 

\begin{example}\label{graded unip group example}
    Let \(\Uhat\subset \operatorname{SL}_2\) be the subgroup of upper triangular matrices. Then \(\Uhat\) is a graded unipotent group. Indeed, \(\Uhat\) is the parabolic subgroup associated to the cocharacter \(\lambda:\mathbb{G}_m\to \operatorname{SL}_2\) defined by 
    \[t\mapsto \begin{pmatrix}
        t & 0\\ 0 & t^{-1}
    \end{pmatrix}\]
    and the associated Levi subgroup is \(\mathbb{G}_m\).
    The semi-direct product decomposition provided by \Cref{semi-dir decomp} is
    \[\Uhat=\mathbb{G}_a\rtimes \mathbb{G}_m\]
    where \(\mathbb{G}_m\) acts on \(\mathbb{G}_a\) with weight \(2\).
\end{example}
\section{Moduli spaces for $\Theta$-strata}
In this section we prove \Cref{Theorem A} after introducing and studying the necessary technical tools.

\subsection{Contracting action on filtrations}
We will use the geometry of \({\Theta}\) to induce a contracting \({\mathbb{G}_m}\) action on the fibers of the morphism \({\gr:\Filt\to \Grad}\). 
We will give two descriptions of this action and prove that they agree in \Cref{action comparison}. We call the first one the rescaling action and the second one the canonical action.

For the rescaling action, the key observation is that multiplication on \({\mathbb{A}^1}\) defines a family of endomorphisms of \({\Theta}\) degenerating the identity \({\Theta\to \Theta}\) to \({\Theta\to B\mathbb{G}_m\to \Theta}\) given by \({(\mathcal{L},s)\mapsto (\mathcal{L}, 0)}\) on points, and that this degeneration is constant when restricted to \({B\mathbb{G}_m\subset \Theta}\). Given a filtration \({f:\Theta\to \mathcal{X}}\), this family induces, by precomposition, a degeneration of \({f}\) to the associated split filtration, \({\sigma\circ \gr(f)}\). 
This description of the action is essentially the morphism used in \cite[Proposition 1.3.8]{halpernleistner2022structure} to prove that \({\Filt}\) admits a deformation retraction onto \({\Grad}\).\footnote{It is the restriction of the \({\Theta}\)-deformation retraction defined in \cite[Proposition 1.3.8]{halpernleistner2022structure} to \({\mathbb{A}^1\times \Filt}\).}

The canonical action comes from the cocharacter \(\lambda:\mathbb{G}_m\to \text{Aut}(B\mathbb{G}_m\to \mathcal{X})\) obtained from describing maps \(B\mathbb{G}_m\to \mathcal{X}\) in terms of descent data. This cocharacter induces a natural action of \(\mathbb{G}_m\) on the fiber of \({\gr:\Filt \to \Grad}\).
We will see in \Cref{gmcan properties} and \Cref{canonical action and canonical subgroup} that this action factors through a central subgroup of the inertia stack, \({\mathcal{I}_{\Grad}}\).


Let 
\[\rs:\mathbb{A}^1\times \Theta\to \Theta\]
be defined on \({T}\)-points for any \({S}\)-scheme \({T}\) by
\[(t,(\mathcal{L},s))\mapsto (\mathcal{L},t\cdot s)\]
that is we multiply \({s:\mathcal{L}\to \mathcal{O}_{T}}\) by \({t\in H^0(T,\mathcal{O}_T)=\mathbb{A}^1(T)}\). We call this morphism \({\rs}\) because it "rescales" the morphism \({s:\mathcal{L}\to \mathcal{O}_T}\). 

This induces the degeneration of the identity morphism \({\id:\Theta\to \Theta}\) to the composition \({\Theta\to B\mathbb{G}_m\to \Theta}\) in the following way.

Given \({t\in \mathbb{A}^1(T)}\) let 
\[\rs(t,-):\Theta_T\to \Theta_T\]
denote the composition \({\rs\circ (t\times \id_{\Theta_T}):\Theta_T\to \Theta_T}\). In particular for any \({t\in \mathbb{G}_m(T)}\), \({\id_{\Theta}\cong\rs(t,-)}\) via the \({2}\)-morphism \({t^{-1}\cdot:(\mathcal{L},s)\to (\mathcal{L},ts)}\) defined by scalar multiplication by \({t^{-1}}\) for all pairs \({(\mathcal{L},s)}\). We also have \({\rs(0,-)=0\circ q}\), as we defined the morphisms \({q:\Theta\to B\mathbb{G}_m}\) by \({q(\mathcal{L},s)=\mathcal{L}\in B\mathbb{G}_m(T)}\) and \({0:B\mathbb{G}\to \Theta}\) by \({0(\mathcal{L})=(\mathcal{L},0)\in \Theta(T)}\).

Let \({m:\mathbb{A}^1\times \mathbb{A}^1\to \mathbb{A}^1}\) denote the multiplication morphism of \({\mathbb{A}^1}\) as in \Cref{gm-actions-on-spaces}. 

Let us show that \({\rs:\mathbb{A}^1\times \Theta\to \Theta}\) behaves like an \({(\mathbb{A}^1,\cdot)}\)-monoid action on \({\Theta}\), which restricts to the trivial action on \({B\mathbb{G}_m\subset \Theta}\).
\begin{lemma}\label{rs is action}
    The morphism \({\rs:\mathbb{A}^1\times\Theta\to \Theta}\) makes the following diagrams commute strictly\footnote{Strict commutativity means that all implicit \({2}\)-morphisms in the diagrams are the identity \({2}\)-morphism.}:
    \begin{enumerate}
        \item Associativity:
        \[\begin{tikzcd}
	{\mathbb{A}^1\times\mathbb{A}^1\times\Theta} & {\mathbb{A}^1\times \Theta} \\
	{\mathbb{A}^1\times \Theta} & \Theta
	\arrow["{\operatorname{id}_{\mathbb{A}^1}\times \rs}", from=1-1, to=1-2]
	\arrow["{m\times\id_{\Theta}}"', from=1-1, to=2-1]
	\arrow["\rs", from=1-2, to=2-2]
	\arrow["\rs"', from=2-1, to=2-2]
\end{tikzcd}\]
        \item Identity:
        \[\begin{tikzcd}
	{\mathbb{A}^1\times\Theta} & \Theta \\
	\Theta
	\arrow["{ \rs}", from=1-1, to=1-2]
	\arrow["{1\times \id_{\Theta}}", from=2-1, to=1-1]
	\arrow["{\id_{\Theta}}"', from=2-1, to=1-2]
\end{tikzcd}\]
        \item The origin is invariant:
        \[\begin{tikzcd}
	{\mathbb{A}^1\times \Theta} & \Theta \\
	{\mathbb{A}^1\times B\mathbb{G}_m} & {B\mathbb{G}_m}
	\arrow["\rs", from=1-1, to=1-2]
	\arrow["{\id_{\mathbb{A}^1\times 0}}", from=2-1, to=1-1]
	\arrow["{p_2}", from=2-1, to=2-2]
	\arrow["0"', from=2-2, to=1-2]
\end{tikzcd}\]
        \item Action over \({B\mathbb{G}_m}\):
        \[\begin{tikzcd}
	{\mathbb{A}^1\times \Theta} & \Theta \\
	{\mathbb{A}^1\times B\mathbb{G}_m} & {B\mathbb{G}_m}
	\arrow["\rs", from=1-1, to=1-2]
	\arrow["{\id_{\mathbb{A}^1}\times q}"', from=1-1, to=2-1]
	\arrow["q", from=1-2, to=2-2]
	\arrow["{p_2}"', from=2-1, to=2-2]
\end{tikzcd}\]
    \end{enumerate}
\end{lemma}
\begin{proof}
        Let \({T}\) be an \({S}\)-scheme, and \({(\mathcal{L},s)\in \Theta(T)}\).
        The first two properties follows directly from the fact that the sheaf of morphisms \({s:\mathcal{L}\to \mathcal{O}_T}\) is an \({\mathcal{O}_T}\) module, \({\mathcal{L}^{\lor}}\), and that \({\rs(t,(\mathcal{L},s))=(\mathcal{L},ts)}\) where \({ts}\) is defined by the \({\mathcal{O}_T}\)-module structure of \({\mathcal{L}^{\lor}}\).

        Similarly, \({\rs(t,(\mathcal{L},0))=(\mathcal{L},0)}\) for any \({t\in \mathbb{A}^1(T)}\), which implies part 3, and
        
        \({q\circ \rs(t,(\mathcal{L},s))=q(\mathcal{L},st)=\mathcal{L}}\), proving part 4.
\end{proof}

Using \({\rs:\mathbb{A}^1\times \Theta\to \Theta}\), we define 
\[\rsInduced:\mathbb{A}^1\times \Filt\to \Filt\]
by
\[(t,f:\Theta_T\to \mathcal{X})\mapsto f\circ \rs(t,-):\Theta_T\to \mathcal{X}.\]

Let us use \Cref{rs is action} to show that the morphism \({\rsInduced:\mathbb{A}^1\times\Filt\to \Filt}\) also behaves like an \({(\mathbb{A}^1,\cdot)}\)-monoid action.

\begin{lemma}\label{strict-contracting action on Filt over Grad}
    The morphism \({\rsInduced:\mathbb{A}^1\times\Filt\to \Filt}\)
    makes the following diagrams strictly \({2}\)-commutative:
    \begin{enumerate}
        \item Associativity:
        \[\begin{tikzcd}
	{\mathbb{A}^1\times \mathbb{A}^1\times \Filt} && {\mathbb{A}^1\times \Filt} \\
	{\mathbb{A}^1\times \Filt} && \Filt
	\arrow["{\id_{\mathbb{A}^1}\times \rsInduced}", from=1-1, to=1-3]
	\arrow["{m\times \id_{\Filt}}"', from=1-1, to=2-1]
	\arrow["\rsInduced", from=1-3, to=2-3]
	\arrow["\rsInduced"', from=2-1, to=2-3],
\end{tikzcd}\]
        \item Identity: 
        \[\begin{tikzcd}
	{\mathbb{A}^1\times \Filt} && \Filt \\
	\Filt
	\arrow["\rsInduced", from=1-1, to=1-3]
	\arrow["{1\times \id_{\Filt}}", from=2-1, to=1-1]
	\arrow["{\id_{\Filt}}"', from=2-1, to=1-3]
\end{tikzcd}\],
        \item Compatibility with \({\gr:\Filt\to \Grad}\):
        \[\begin{tikzcd}
	{\mathbb{A}^1\times \Filt} && \Filt \\
	{\mathbb{A}^1\times \Grad} && \Grad
	\arrow["\rsInduced", from=1-1, to=1-3]
	\arrow["{\id_{\mathbb{A}^1}\times\gr}"', from=1-1, to=2-1]
	\arrow["\gr", from=1-3, to=2-3]
	\arrow["{p_2}"', from=2-1, to=2-3]
\end{tikzcd}\],
        \item Invariant section:
        \[\begin{tikzcd}
	{\mathbb{A}^1\times \Filt} && \Filt \\
	{\mathbb{A}^1\times \Grad} && \Grad
	\arrow["\rsInduced", from=1-1, to=1-3]
	\arrow["{\id_{\mathbb{A}^1}\times \sigma}", from=2-1, to=1-1]
	\arrow["{p_2}"', from=2-1, to=2-3]
	\arrow["\sigma"', from=2-3, to=1-3]
\end{tikzcd}\].
    \end{enumerate}
    In particular, if \({f_0:T\to \Grad}\) is a morphism such that the fiber of \({\gr:\Filt\to \Grad}\) over \({f_0:T\to \Grad}\) is an algebraic space over \({T}\), then \({\rs}\) defines a contracting \({\mathbb{G}_m}\) action on the fiber and \({\sigma}\) defines an invariant section.
\end{lemma}
\begin{proof}
    Parts 1 to 4 are all formal consequences of parts 1 to 4 of \Cref{rs is action}. To illustrate this, we show part 1. 
    Let \({(t,s,f:\Theta_T\to \mathcal{X})\in\mathbb{A}^1\times \mathbb{A}^1\times \Filt}\), then by definition
    \[\rsInduced\circ (\id_{\mathbb{A}^1}\times \rsInduced)(t,s,f)=f\circ \rs(s,-) \circ(\id_{\mathbb{A}^1}\times \rs(t,-))\]
    which equals
    \[f\circ \rs(st,-)=f\circ \rs(ts,-)\]
    by part 1 of \Cref{rs is action}, and the commutativity of multiplication on \({\mathbb{A}^1}\).
    On the other hand, 
    \[\rsInduced\circ (m\times \id_{\Filt})((t,s,f:\Theta_T\to \mathcal{X}))=f\circ \rs(ts,-)\]
    showing part 1. 
    
    Given a morphism \({f_0:T\to \Grad}\) such that \({\Filt\times_{\Grad}T\to T}\) is an algebraic space, part 3 ensures that the pullback of \({\rsInduced}\),\[{\rsInducedOverFnull:\mathbb{A}^1\times \Filt\times_{\Grad}T\to T}\] is well defined. By part 1 and 2 this induced morphism is a \({(\mathbb{A}^1,\cdot)}\)-monoid action. Hence the restriction \[{\rsInducedOverFnull:\mathbb{G}_m\times \Filt\times_{\Grad}T\to T}\] is a contracting \({\mathbb{G}_m}\)-action by part 5 of \Cref{multip-act}. Finally, part 4 implies that the pullback of \({\sigma:\Grad\to \Filt}\) is an invariant section for this action.
\end{proof}
Next we will give a description of the fixed points of the contracting action 
\[{\rsInducedOverFnull:\mathbb{A}^1\times \FiltZ\times_{\mathfrak{Z}}T\to \FiltZ\times_{\mathfrak{Z}}T}\]
in terms of the section \({\sigma_T}\). We do this by first identifying the essential image of \({\sigma:\Grad\to \Filt}\) with the essential image of \({\rsInduced|_{\{0\}\times \Filt}:\Filt\to \Filt}\) where we identify the stacks with their underlying fibered categories.

\begin{lemma}\label{essential fixed points}\label{representable fixed points}
    The essential image of \({\rsInduced|_{\{0\}\times \Filt}:\Filt\to \Filt}\) and the essential image of \({\sigma:\Grad\to \Filt}\) are the same subcategory of the underlying fibered category of the stack \({\Filt}\).

    In particular, if \({f_0:T\to \Grad}\) is a morphism such that the fiber of \({{\gr:\Filt\to \Grad}}\) over \({f_0:T\to \Grad}\) is an algebraic space over \({T}\), then the section 
    \[{{\sigma_T:T\to \Filt\times_{\Grad} T}}\] 
    identifies \({T}\) with the fixed point subspace of the contracting \({\mathbb{G}_m}\)-action on the fiber.
\end{lemma}
\begin{proof}
    Let \({T}\) be a \({S}\)-scheme, we will prove that the essential images of \({\sigma:\Grad\to \Filt}\) and \({\rs(0,-):\Filt\to \Filt}\) are equal on \({T}\)-points.
    Let \({g\in \Grad(T)}\), by definition \({\sigma(g)=g\circ q}\), where \({q:\Theta\to B\mathbb{G}_m}\) is defined by \({q(\mathcal{L},s)=\mathcal{L}}\) for any pair \({(\mathcal{L},s)\in \Theta(T)}\). 
    Let \({f\in \Filt(T)}\), then \({\rsInduced(0,f)=f\circ 0\circ  q}\) where \({0:B\mathbb{G}_m\to \Theta}\) is defined by \({0(\mathcal{L})=(\mathcal{L},0)}\) for any line bundle \({\mathcal{L}\in B\mathbb{G}_m(T)}\), in particular, \({\rsInduced(0,f)=\sigma(f\circ 0)}\) and \({\sigma(g)=\rsInduced(0,\sigma(g))}\) since \({q\circ 0=\id_{B\mathbb{G}_m}}\)
    so the essential image of \({\rs(0,-)}\) is equal to the essential image of \({\sigma}\). 

    Let \({f_0:T\to \Grad}\) be a morphism such that the fiber of \({{\gr:\Filt\to \Grad}}\) is an algebraic space over \({T}\).
    Since \({\Filt\times_{\Grad}T}\) is an algebraic space and \({\sigma_T:T\to \Filt\times_{\Grad}T}\) is a section, it is fully faithful as a functor between categories fibered in setoids (groupoids equivalent to sets). Hence \({T}\) is isomorphic via \({\sigma_T}\) to its essential image, which is equal to the essential image of \({\rsInducedOverFnull(0,-)}\) by the first part of this lemma. 
    Since the induced action of \({\mathbb{G}_m}\) on \({\Filt\times_{\Grad}T}\) is contracting, the morphism \({\plim:\Filt\times_{\Grad}T\to (\Filt\times_{\Grad}T)^0}\), is a split surjection onto the fixed points and \({\plim=\rsInducedOverFnull(0,-)}\). 
    In particular, \({\sigma_T:T\to (\Filt\times_{\Grad}T)^0}\) is an isomorphism with inverse \({\gr_T|_{(\Filt\times_{\Grad}T)^0}:(\Filt\times_{\Grad}T)^0\to T}\).
\end{proof}

We now turn to the definition of the canonical \({\mathbb{G}_m}\)-action on the fibers and its identification with the rescaling action.

Let \({T}\) be a \({S}\)-scheme and let \({f_0:T\to \Grad}\) be defined by a morphism \({B\mathbb{G}_{m,T}\to \mathcal{X}}\). To give such a morphism is equivalent to giving a \({T}\)-point \({x_{f_0}:T\to \mathcal{X}}\) and a cocharacter \({\lambda_{f_0}:\mathbb{G}_{m,T}\to \text{Aut}_{\mathcal{X}}(x_{f_0})}\). The correspondence goes via the atlas \({T\to B\mathbb{G}_{m,T}}\) and descent. Étale locally on \(T\) the morphism \(f_0:B\mathbb{G}_m\to \mathcal{X}\) induced by \((x_{f_0},\lambda_{f_0})\) is described by the morphism of groupoids induced by the commutative diagram
\[\begin{tikzcd}
	& {\text{Aut}_{\mathcal{X}}(x_{f_0})} \\
	{\mathbb{G}_{m,T}} & {X\times_{\mathcal{X}}X} \\
	T & X \\
	& {\mathcal{X}}
	\arrow[hook, from=1-2, to=2-2]
	\arrow["{\lambda_{f_0}}", from=2-1, to=1-2]
	\arrow[from=2-1, to=2-2]
	\arrow[shift right, from=2-1, to=3-1]
	\arrow[shift left, from=2-1, to=3-1]
	\arrow[shift right, from=2-2, to=3-2]
	\arrow[shift left, from=2-2, to=3-2]
	\arrow[from=3-1, to=3-2]
	\arrow["{x_{f_0}}"', from=3-1, to=4-2]
	\arrow[from=3-2, to=4-2]
\end{tikzcd}\]
where \(X\to \mathcal{X}\) is a smooth atlas.

Note that \({\text{Aut}_{\Grad}(f_0)}\) is naturally identified with the centralizer of \({\lambda_{f_0}}\) in \({\text{Aut}_{\mathcal{X}}(x_{f_0})}\) as a group scheme over \({T}\). Indeed an automorphism of \({f_0:B\mathbb{G}_{m,T}\to \mathcal{X}}\) restricts to an automorphism of \({x_{f_0}:T\to \mathcal{X}}\) that respects the descent datum descending \({x_{f_0}}\) to \({f_0}\), that is, it commutes with \({\lambda_{f_0}:\mathbb{G}_{m,T}\to \text{Aut}_{\mathcal{X}}(x_{f_0})}\).

Using the above, we define
\[\mu^{\operatorname{can}}:\mathbb{G}_m\times \Filt_T \to \Filt_T\]
by
\[(t,(f,\alpha:\gr(f)\cong f_0))\mapsto (f,\alpha \circ  \lambda_{\gr(f)}(t^{-1}):\gr(f)\cong f_0)\]
on objects \({(t,(f,\alpha:\gr(f)\cong f_0))\in\mathbb{G}_m\times \Filt_T}\) over any scheme defined over \({T}\).
If \({\Filt_T}\) is an algebraic space, then this is naturally a \({\mathbb{G}_m}\) action on \({\Filt_T}\) by the basic properties of composition of \({2}\)-morphisms.

\begin{lemma}\label{action comparison}
    Let \({\mathcal{X}}\) be an algebraic stack over \({S}\) and let \({f_0:T\to \Grad}\) be a graded \({T}\)-point of \({\mathcal{X}}\), where \({T}\) is a \({S}\)-scheme. Then there is an isomorphism
    \[\mu^{\operatorname{can}}\cong \rsInducedOverFnull\]
    in the groupoid \({\text{Hom}(\mathbb{G}_m\times \Filt_T,\Filt_T)}\).
\end{lemma}
\begin{proof}
    We need to construct an invertible natural transformation \({{\eta:\mu_{f_0}^{\operatorname{can}}\to \rsInducedOverFnull}}\).
    Given an object \({(t,(f,\alpha))\in \mathbb{G}_m\times \Filt_T}\) where \({\alpha:\gr(f)\cong f_0}\) is an isomorphism in \({\Grad}\), recall that
    \[\mu^{\operatorname{can}}(t,(f,\alpha))=(f,\alpha\circ \lambda_{\gr(f)}(t^{-1}))\]
    and 
    \[\rsInducedOverFnull(t,(f,\alpha))=(\rsInduced(t,f),\alpha)\]
    so we want a morphism
    \[(f,\alpha\circ \lambda_{\gr(f)}(t^{-1}))\to(\rsInduced(t,f),\alpha)\]
    in \({\Filt_T}\).  
    
    Recall that for any \(t\in \mathbb{G}_m(T)\), there is a morphism
    \[t^{-1}\cdot:\id_{\Theta_T}\to \rs(t,-)\]
    in \({\text{Hom}(\Theta_T,\Theta_T)}\)
    which when restricted to \({\text{Hom}(B\mathbb{G}_{m,T},B\mathbb{G}_{m,T})}\) is an automorphism of the identity.
    Consider the diagram
\[\begin{tikzcd}
	{\Theta_T} && {\Theta_T} && {\mathcal{X}} \\
	\\
	{B\mathbb{G}_{m,T}} && {B\mathbb{G}_{m,T}}
	\arrow[""{name=0, anchor=center, inner sep=0}, "\id", curve={height=-12pt}, from=1-1, to=1-3]
	\arrow[""{name=1, anchor=center, inner sep=0}, "{\rs(t,-)}"', curve={height=12pt}, from=1-1, to=1-3]
	\arrow[""{name=2, anchor=center, inner sep=0}, from=1-3, to=1-5]
	\arrow["0", from=3-1, to=1-1]
	\arrow[""{name=3, anchor=center, inner sep=0}, "\id"', curve={height=12pt}, from=3-1, to=3-3]
	\arrow[""{name=4, anchor=center, inner sep=0}, "\id", curve={height=-12pt}, from=3-1, to=3-3]
	\arrow["0", from=3-3, to=1-3]
	\arrow[""{name=5, anchor=center, inner sep=0}, "{f_{0}}"', from=3-3, to=1-5]
	\arrow["t^{-1}\cdot", shorten <=3pt, shorten >=3pt, Rightarrow, from=0, to=1]
	\arrow["\alpha", shorten <=4pt, shorten >=4pt, Rightarrow, from=2, to=5]
	\arrow["\text{id}", shorten <=3pt, shorten >=3pt, Rightarrow, from=4, to=3]
\end{tikzcd}\]
which \(2\)-commutes by the basic properties of \({\rs(t,-)}\) outlined in the beginning of this section. 

This induces a commutative diagram 
\[\begin{tikzcd}
	\Theta && {\mathcal{X}} \\
	\\
	{B\mathbb{G}_m}
	\arrow[""{name=0, anchor=center, inner sep=0}, "f"{description}, curve={height=-12pt}, from=1-1, to=1-3]
	\arrow[""{name=1, anchor=center, inner sep=0}, "{\rs(t,f)}"{description}, curve={height=12pt}, from=1-1, to=1-3]
	\arrow[from=3-1, to=1-1]
	\arrow[""{name=2, anchor=center, inner sep=0}, "{f_0}"', from=3-1, to=1-3]
	\arrow["{t^{-1}\cdot}"', shorten <=3pt, shorten >=3pt, Rightarrow, from=0, to=1]
	\arrow["\alpha"', shorten <=3pt, shorten >=3pt, Rightarrow, from=1, to=2]
\end{tikzcd}\]
i.e we get a morphism \(t^{-1}\cdot:f\cong\rs_*(t,f)\) such that 
\[\begin{tikzcd}
	{\gr(f)} & {\gr(\rs_*(t,f))} \\
	& {f_0}
	\arrow["{\gr(t^{-1}\cdot)}", from=1-1, to=1-2]
	\arrow["{\alpha\circ \lambda_{\gr(f)}(t^{-1})}"', from=1-1, to=2-2]
	\arrow["\alpha", from=1-2, to=2-2]
\end{tikzcd}\]
commutes. By the definition of fiber products of stacks, see e.g \cite[Tag 003O]{stacks-project}, we get a morphism
\[\eta_{(t,(f,\alpha))}:\mu^{\operatorname{can}}(t,(f,\alpha))\to \rsInducedOverFnull(t,(f,\alpha))\]
in \(\Filt_T\).

To finish the proof, we need to show that this collection of morphisms is natural in \((t,(f,\alpha))\in \mathbb{G}_m\times \Filt_T\). As \({\mathbb{G}_m}\) and \({T}\) are schemes, all morphisms in \({\mathbb{G}_m\times\Filt_T}\) are of the form 
\[\beta:(t,(g,\alpha\circ \gr(\beta)))\to (t,(f,\alpha))\] 
where \({\beta:g\to f}\) is a morphism in \({\Filt}\), \cite[Tag 02X9]{stacks-project}. 
As \({\eta_{(t,(f,\alpha))}}\) is defined by applying \({f}\) to the isomorphism \({t^{-1}\cdot:\id_{\Theta}\to \rs(t,-)}\), 
to show that \({\eta}\) is a natural transformation, it is enough to check that the diagram
\[
    \begin{tikzcd}
	{g} & {\rsInduced(t,g)} \\
	{f} & {\rsInduced(t,f)}
	\arrow["f(t^{-1}\cdot)"', from=2-1, to=2-2]
	\arrow["{\rsInduced(t,\beta)}", from=1-2, to=2-2]
	\arrow["\beta"', from=1-1, to=2-1]
	\arrow["g(t^{-1}\cdot)", from=1-1, to=1-2]
 \end{tikzcd}
\]
commutes. By definition of \({\rsInduced}:\mathbb{A}^1\times \Filt\to \Filt \), this is the same as showing that for any object \({x\in \Theta_T}\), the diagram,
\[    \begin{tikzcd}
	{g(x)} & {g\circ \rs(t,x)} \\
	{f(x)} & {f\circ \rs(t,x)}
	\arrow["f(t^{-1}\cdot)"', from=2-1, to=2-2]
	\arrow["{\beta_{\rs(t,x)}}", from=1-2, to=2-2]
	\arrow["\beta_x"', from=1-1, to=2-1]
	\arrow["g(t^{-1}\cdot)", from=1-1, to=1-2]
 \end{tikzcd}\]
  commutes, which it does as \({\beta:g\to f}\) is a natural transformation.
\end{proof}

\subsection{Descending the action along rigidifications}

In this section we show that the action on the fibers of \({\gr:\Filt\to \Grad}\) defined in the previous section descends to an action on the fibers of \({\FiltZRig \to \mathfrak{Z}\subset \Grad}\), whenever the rigidification \({\rho:\FiltZ\to\FiltZRig }\) is defined. 

To rigidify with respect to the relative inertia of \({\gr:\Filt\to \Grad}\), \({\mathcal{I}_{\gr}}\), we want to apply \cite[Appendix A]{Tame-stacks}. To do this, we need \({\mathcal{I}_{gr}\subset \mathcal{I}_{\FiltZ}}\) to be a closed substack that is flat and of finite presentation over \({\FiltZ}\). 

We first show that in the situation we are interested in, \({\mathcal{I}_{\gr}\subset \mathcal{I}_{\mathcal{X}}}\) is always a closed substack.

\begin{lemma}\label{relative inertia is closed}
    Let \({f:\mathcal{X}\to \mathcal{Y}}\) be a morphism of algebraic stacks. If \({\mathcal{Y}}\) is has a separated diagonal, then \({\mathcal{I}_{f}\to \mathcal{I}_{\mathcal{X}}}\) is a closed immersion.
\end{lemma}
\begin{proof}
    The following diagram is cartesian
\[\begin{tikzcd}
	{\mathcal{I}_{f}} & {\mathcal{I}_{\mathcal{X}}} \\
	\mathcal{X} & {\mathcal{X}\times_{\mathcal{Y}}\mathcal{I}_{\mathcal{Y}}}
	\arrow[from=1-1, to=1-2]
	\arrow[from=1-1, to=2-1]
	\arrow[from=1-2, to=2-2]
	\arrow[from=2-1, to=2-2]
\end{tikzcd}\]
    where \({\mathcal{X}\times_{\mathcal{Y}}\mathcal{I}_{\mathcal{Y}}}\) is the section of the first projection defined by \({x\mapsto (x,(f(x),\id_{f(x)}))}\), and\\ \({\mathcal{I}_{\mathcal{X}}\to \mathcal{X}\times_{\mathcal{Y}}\mathcal{I}_{\mathcal{Y}}}\) is defined by \({(x,\alpha)\mapsto (x,(f(x),f(\alpha)))}\). Hence it is enough to show that the section \({\mathcal{X}\to \mathcal{X}\times_{\mathcal{Y}}\mathcal{I}_{\mathcal{Y}}}\) is a closed immersion, which holds if the projection \({\mathcal{X}\times_{\mathcal{Y}}\mathcal{I}_{\mathcal{Y}}\to \mathcal{X}}\) is separated, and in particular if \({\mathcal{I}_{\mathcal{Y}}\to \mathcal{Y}}\) is separated.
    
    By \cite[Tag 0CL0]{stacks-project}, \({\mathcal{I}_{\mathcal{Y}}\to \mathcal{Y}}\) is separated if and only \({\mathcal{Y}}\) has separated diagonal, which is true by assumption and we are done.
\end{proof}

By \Cref{basic-prop-of-filt-grad}, whenever \({\mathcal{X}}\) is an algebraic stack with separated and quasi-compact diagonal locally of finite presentation over \({S}\), the same holds for \({\Grad}\).
Thus, applying \Cref{relative inertia is closed}, whenever \({\mathcal{I}_{\gr}}\) is flat and of finite presentation over \({\FiltZ}\), the rigidification \({\FiltZRig}\) is an algebraic stack equipped with a morphism \({\rho:\FiltZ\to \FiltZRig}\) making \({\FiltZ}\) a fppf-gerbe over \({\FiltZRig}\) by \cite[Theorem A.1]{Tame-stacks}.

\begin{lemma}\label{action respects inertia}
     Let \({\mathcal{X}}\) be an algebraic stack locally of finite presentation, with separated and quasi-compact diagonal and affine stabilizers over a quasi-separated scheme \({S}\).
    If \({{\mathcal{I}_{\gr}\to \FiltZ}}\) is flat of finite presentation, then the morphism \({{\rsInduced:\mathbb{A}^1\times \FiltZ\to \FiltZ}}\) descends to a morphism \[\rsInducedRig:\mathbb{A}^1\times \FiltZRig\to \FiltZRig.\]
\end{lemma}
\begin{proof}
    By the universal property of rigidification, namely a morphism \({\mathcal{X}\to \mathcal{Y}}\) descends to\\ \({\mathcal{X}\!\!\!\fatslash H \to \mathcal{Y}}\) if and only if \({\mathcal{I}_{\mathcal{X}/\mathcal{Y}}\subset H}\) (\Cref{universal property of rigi}), it is enough to show that the map on inertia induced by \({\rsInduced:\mathbb{A}^1\times \FiltZ\to \FiltZ}\) induces a morphism \({\mathbb{A}^1\times \mathcal{I}_{\gr}\to\mathcal{I}_{\gr}}\).
    
    A \({T}\)-point of \({\mathbb{A}^1\times \mathcal{I}_{\text{gr}}}\) is a pair \({(t,\begin{tikzcd}
	{\Theta_T}
	\arrow[bend left=20,"f", shift left]{r}[name=U]{}
	\arrow[bend right=20,"f"', shift right]{r}[name=L]{} 
    \arrow[ Rightarrow,"\alpha",to path=(L) -- (U)\tikztonodes]{}
 & {\mathcal{X}}
\end{tikzcd})}\) where \({t\in \mathbb{A}^1(T)}\), \({f\in \FiltZ(T)}\) and\\ \({\alpha\in \text{Aut}_{\FiltZ}(f)}\). Applying \({\gr\circ \rsInduced:\mathbb{A}^1\times \FiltZ\to \mathfrak{Z}}\) to this pair yields the diagram 
\[\begin{tikzcd}
    B\mathbb{G}_{m,T} \arrow["0"]{r}{}
    &
    \Theta_T\arrow["t\times \text{id}_{\Theta}"]{r}{}
    &
    \mathbb{A}^{1}\times\Theta_T \arrow["\rs"]{r}{}
    &
	{\Theta_T}
	\arrow[bend left=20,"f", shift left]{r}[name=U]{}
	\arrow[bend right=20,"f"', shift right]{r}[name=L]{} 
    \arrow[ Rightarrow,"\alpha",to path=(L) -- (U)\tikztonodes]{}
 & {\mathcal{X}}
\end{tikzcd}\]
describing a point of \({\mathcal{I}_{\mathfrak{Z}}}\). As the composition \({\rs \circ (t\times \text{id}_{\Theta})\circ 0=\text{id}_{B\mathbb{G}_m}}\), the point is equal to the one described by 
\[\begin{tikzcd}
    B\mathbb{G}_{m,T} \arrow["0"]{r}{}
    &
	{\Theta_T}
	\arrow[bend left=20,"f", shift left]{r}[name=U]{}
	\arrow[bend right=20,"f"', shift right]{r}[name=L]{} 
    \arrow[ Rightarrow,"\alpha",to path=(L) -- (U)\tikztonodes]{}
 & {\mathcal{X}}
\end{tikzcd}\]
which is equal to the one described by applying \({\text{gr}:\FiltZ\to \mathfrak{Z}}\) to
\[\begin{tikzcd}
	{\Theta_T}
	\arrow[bend left=20,"f", shift left]{r}[name=U]{}
	\arrow[bend right=20,"f"', shift right]{r}[name=L]{} 
    \arrow[ Rightarrow,"\alpha",to path=(L) -- (U)\tikztonodes]{}
 & {\mathcal{X}}
\end{tikzcd}.\]
But this is a point in \({\mathcal{I}_{\text{gr}}}\) by assumption, so \({\alpha}\) induces the identity in \({\text{Aut}_{\mathfrak{Z}}(\text{gr}(f))}\), in other words \({\rsInduced(t,\alpha)\in \mathcal{I}_{\gr}(T)}\). 
\end{proof}
Let \({\grRig:\FiltZRig\to \mathfrak{Z}}\) be the morphism induced by the universal property of rigidification applied to \({\gr:\FiltZ\to \mathfrak{Z}}\), i.e. it is the unique morphism \({\FiltZRig\to \mathfrak{Z}}\) such that \({\gr=\grRig\circ \rho}\).

\begin{lemma}\label{action on fibers 2}
     If \({\mathcal{I}_{\text{gr}}}\) is flat and finitely presented over \({\FiltZ}\), then 
     \[{\rsInducedRig:\mathbb{A}^1\times \FiltZRig\to \FiltZRig}\] 
     induces a contracting \({\mathbb{G}_m}\) action on the fibers of \({\grRig:\FiltZRig\to \mathfrak{Z}}\), and the section\\ \({\rho\circ \sigma:\mathfrak{Z}\to \FiltZRig}\) is invariant for this action.
\end{lemma}
\begin{proof}
    Let \({f_0:T\to \mathfrak{Z}}\) be a morphism from a scheme \({T}\).
    As was shown in the proof of \Cref{strict-contracting action on Filt over Grad}, the morphism 
    \[\rsInducedOverFnull:\mathbb{A}^1\times \FiltZ\times_{\mathfrak{Z}}T\to \FiltZ\times_{\mathfrak{Z}}T\]
    satisfies the axioms of a monoid action, that is associativity and identity. 
    To see that this is also true for 
    \[\rsInducedFnullRig:\mathbb{A}^1\times \FiltZRig\times_{\mathfrak{Z}}T\to \FiltZRig\times_{\mathfrak{Z}}T\]
    we apply \Cref{formal corollary} to the diagrams describing associativity and identity of
    \({\rsInducedOverFnull}\).
    To do this, we rigidify \({\mathbb{A}^1\times \FiltZ\times_{\mathfrak{Z}}T}\) with respect to \({\mathbb{A}^1\times \mathcal{I}_{\text{gr}}\times_{\mathfrak{Z}}T}\), \({\mathbb{A}^1\times \mathbb{A}^1\times \FiltZ\times_{\mathfrak{Z}}T}\) with respect to \({\mathbb{A}^1\times \mathbb{A}^1\times \mathcal{I}_{\text{gr}}\times_{\mathfrak{Z}}T}\) and \({\FiltZ\times_{\mathfrak{Z}}T}\) with respect to \({\mathcal{I}_{\text{gr}}\times_{\mathfrak{Z}}T}\).
    
The invariance of \({\rho \circ \sigma:\mathfrak{Z}\to \FiltZRig}\) follows from the invariance of \({\sigma}\). Indeed, for all \({t\in \mathbb{A}^1(T)}\), \({\rsInduced\circ (t\times \sigma)=\sigma}\), and by the universal property of rigidification (\Cref{universal property of rigi}),\\ 
\({\rho\circ \rsInduced=\rsInducedRig\circ (\text{id}\times \rho)}\). Hence for all \({t\in \mathbb{A}^1(T)}\), the diagram 
\[\begin{tikzcd}
	{\mathfrak{Z}} && {\mathbb{A}^1\times \FiltZRig} \\
	&& \FiltZRig
	\arrow["{t\times (\rho\circ \sigma)}", from=1-1, to=1-3]
	\arrow["{\rho\circ\sigma}"', from=1-1, to=2-3]
	\arrow["\rsInducedRig", from=1-3, to=2-3]
\end{tikzcd}\]
commutes strictly, finishing the proof.
\end{proof}

\begin{lemma}\label{fixed points in FiltRig}
    Let \({\mathfrak{Z}\subset \Grad}\) be a substack such that \({\gr:\FiltZ\to \mathfrak{Z}}\) has relative inertia flat and of finite presentation over \({\FiltZ}\).
    Then, given any morphism \({f_0:T\to \mathfrak{Z}}\) from a scheme \({T}\), the section 
    \[\rho\circ \sigma_T:T\to \left(\FiltZRig\right)\times_{\mathfrak{Z}}T\]
    of \({\grRig:\FiltZRig\to \mathfrak{Z}}\) pulled back along \({T\to \mathfrak{Z}}\) identifies \({T}\) with the fixed point subspace \({\left(\left(\FiltZRig\right)\times_{\mathfrak{Z}}T\right)^0\subset \left(\FiltZRig\right)\times_{\mathfrak{Z}}T}\).
\end{lemma}
\begin{proof}  
    As \({\left(\FiltZRig\right)\times_{\mathfrak{Z}}T}\) is an algebraic space with a contracting \({\mathbb{G}_m}\)-action, the fixed points are given by the image of \({\rsInducedRig(0,-)}\).
    The essential image of 
    \[\rsInducedFnullRig(0,-):\left(\FiltZRig\right)\times_{\mathfrak{Z}}T\to \left(\FiltZRig\right)\times_{\mathfrak{Z}}T\] 
    is the essential image of 
    \[\rho \circ \rsInducedOverFnull(0,-):\FiltZ\times_{\mathfrak{Z}}T\to \left(\FiltZRig\right)\times_{\mathfrak{Z}}T.\] 
    By \Cref{essential fixed points}, this is equal to the essential image of 
    \[\rho\circ \sigma:T\to \left(\FiltZRig\right)\times_{\mathfrak{Z}}T\]
    and as \({\rho\circ \sigma}\) is a section, it is fully faithfull as a functor hence an equivalence of categories onto the essential image of \({\rsInducedFnullRig}\), giving the identification.
\end{proof}

\subsection{The rigidified center}

In this section we study a natural subgroup of the inertia group stack of \({\Grad}\), which arises from the universal morphism \({B\mathbb{G}_m\to \mathcal{X}}\) over \({\Grad}\). 
We use this subgroup stack to define the rigidified center \({\Zrig}\) of the stack of filtrations centered at \({\mathfrak{Z}}\), and relate this subgroup to the contracting action on the fibers of \({\gr:\Filt\to \Grad}\). 

We will define the subgroup stack of the inertia stack of \({\Grad}\), \({\GmCan \subset \mathcal{I}_{\Grad}}\) via a cocharacter of the inertia defined over \({\Grad}\), the universal version of the cocharacter used to define the canonical action above. 

Let
\[\lambda^{\operatorname{can}}:\mathbb{G}_m\times\Grad\to \mathcal{I}_{\Grad}\]
be defined by 
\[(t,g)\mapsto (g,g(t))\]
for \({(t,g)\in \mathbb{G}_m\times \Grad(T)}\) for a \({S}\)-scheme \({T}\), where \({g(t)}\) denotes \({g}\) applied to the automorphism \({t\cdot:\id_{B\mathbb{G}_{m,T}}\to\id_{B\mathbb{G}_{m,T}}}\), which is an automorphism of the morphism  \({B\mathbb{G}_{m,T}\to\mathcal{X}}\) defining \({g}\).
Note that this cocharacter factors through the center of \({\mathcal{I}_{\Grad}}\) as the automorphism group of \({g:B\mathbb{G}_{m,T}\to\mathcal{X}}\) is the centralizer of the cocharacter \({\lambda_g:\mathbb{G}_{m,T}\to \text{Aut}_{\mathcal{X}}(g\circ p)}\) where \({p:T\to B\mathbb{G}_{m,T}}\) is the atlas and \({\lambda(t)= g(t)}\).

We can now prove the first part of \Cref{Theorem A}.

\begin{proposition}\label{gmcan properties}
    Let \({\mathcal{X}}\) be an algebraic stack locally of finite presentation, with separated and quasi-compact diagonal and affine stabilizers over a quasi-separated scheme \({S}\).\\
    Then \({\lambda^{\operatorname{can}}:\mathbb{G}_m\times \Grad\to \mathcal{I}_{\Grad}}\) factors through a closed subgroup stack of \({\mathcal{I}_{\Grad}}\) which we call \({\GmCan}\). This subgroup stack is flat and of finite presentation over \({\Grad}\).
    
    If \({f_0:T\to \Grad}\) is a morphism from a connected scheme, then the pullback of \({\GmCan}\), is either isomorphic to \({\mathbb{G}_{m,T}}\) or the trivial group. The pullback is trivial if and only if the morphism \({B\mathbb{G}_{m,T}\to \mathcal{X}}\) defining \({f_0}\) factors through \({T}\).
\end{proposition}
\begin{proof}
    By part 1 of \Cref{basic-prop-of-filt-grad}, \({\Grad}\) is an algebraic stack with separated and quasi-compact diagonal, with affine stabilizers and locally of finite presentation over  \({S}\). In particular \({\mathcal{I}_{\Grad}\to \Grad}\) is separated.
    Let \({T\to \Grad}\) be an atlas, then the pullback of 
    \[\lambda^{\operatorname{can}}:\mathbb{G}_m\times \Grad\to \mathcal{I}_{\Grad}\]
    to \({T}\) is a morphism of separated group algebraic spaces over \({T}\). By \cite[Proposition 1.3.9]{halpernleistner2022structure}, the kernel of \({\lambda^{\operatorname{can}}_T}\) is locally constant on \({T}\), and thus isomorphic to \({\mu_n}\) for some \({n}\) or \({\mathbb{G}_m}\) on each component. In particular the quotient \({\GmCan_{T}:=\mathbb{G}_{m,T}/\ker(\lambda^{\operatorname{can}})}\) is a flat group scheme over \({T}\), and the induced morphism \({\GmCan_{T}\to \mathcal{I}_{\Grad}\times_{\Grad}T}\) is a monomorphism. Furthermore on each connected component of \({T}\), \({\GmCan_T}\) is either isomorphic to \({\mathbb{G}_m}\) if \({\ker(\lambda^{\operatorname{can}}_T)=\mu_n}\) for some \({n}\) or it is trivial if the kernel is \({\mathbb{G}_{m,T}}\).
    As \({\mathcal{I}_{\Grad}\times_{\Grad}T}\) is separated and \({\GmCan_{T}\subset \mathcal{I}_{\Grad}\times_{\Grad}T}\) is central and in particular normal, it follows that it is a closed subgroup algebraic space, hence \({\GmCan_T}\) descends to a closed subgroup stack of \({\mathcal{I}_{\Grad}}\) with the desired properties.
\end{proof}

\begin{lemma}\label{canonical action and canonical subgroup}
    Let \({f_0:T\to \Grad}\) be a morphism from a \({S}\)-scheme \({T}\).\\
    The action \({\mu^{\operatorname{can}}:\mathbb{G}_{m,T}\times \Filt_T\to \Filt_T}\) defined above \Cref{action comparison} factors through\\ 
    \({{\GmCan_T\times \Filt_T\to \Filt_T}}\).
     In particular, if \({\GmCan}\) is trivial over \({T}\), then \[{\mu^{\operatorname{can}}:\mathbb{G}_{m,T}\times \Filt_T\to \Filt_T}\] is just the second projection i.e. \({\mu^{\operatorname{can}}\cong\rsInducedOverFnull}\) is the trivial \({\mathbb{G}_m}\)-action.
    \end{lemma}
    \begin{proof}
    Define 
    \[\mathbb{G}_{m,T}\times \Filt_T\to \GmCan_T\times \Filt_T\]
    by 
    \[(t,f,\alpha)\mapsto (\lambda_{\gr(f)}(t),f,\alpha)\]
    and note that 
    \[\mu^{\operatorname{can}}(t,f,\alpha)=(f,\alpha\circ \lambda_{\gr(f)}(t^{-1}))\]
    so we get the desired factorization if we define
    \[\GmCan_T\times \Filt_T\to \Filt_T\]
    by 
    \[(\beta,f,\alpha)\mapsto (f,\alpha\circ \beta^{-1}) \]
    i.e by restricting the natural left action of \({\mathcal{I}_{\Grad}\times_{\Grad}T}\) on \({\Filt_T}\) to \({{\GmCan\times_{\Grad}T}}\).
    \end{proof}

\begin{definition}\label{rigidified center}
    Given a substack \({\mathfrak{Z}\to \Grad}\), the \emph{rigidified center of \(\FiltZ\)} is defined as
\[\Zrig:=\mathfrak{Z}\!\!\!\fatslash \GmCan,\]
the rigidification of \({\mathfrak{Z}}\) with respect to \({\GmCan}\).
\end{definition}
The above is well defined as \Cref{gmcan properties} ensures that the assumptions of \cite[Theorem A.1]{Tame-stacks} are satisfied.

\begin{example}\label{Gm-can/RigCent example}
    Let \({\mathcal{X}=B\mathbb{G}_{m}}\) be the classifying stack of the multiplicative group over some scheme \({S}\), then \({\Grad=\coprod_{n\in \mathbb{Z}} B\mathbb{G}_m}\) by \cite[Theorem 1.4.7]{halpernleistner2022structure}, where the component indexed by \({n}\) corresponds to the cocharacter \({\lambda_n:\mathbb{G}_m\to \mathbb{G}_m}\) defined by \({t\mapsto t^n}\), it follows that \({\mathcal{I}_{\Grad}\cong\mathbb{G}_m\times \left(\coprod_{n\in \mathbb{Z}} B\mathbb{G}_m\right)}\). In this case, 
    \[\lambda^{\operatorname{can}}:\mathbb{G}_m\times \left(\coprod_{n\in \mathbb{Z}} B\mathbb{G}_m\right)\to \mathbb{G}_m\times \left(\coprod_{n\in \mathbb{Z}} B\mathbb{G}_m\right)\] 
    restricts to \({\lambda_{n}\times \id_{B\mathbb{G}_m}}\) on the component indexed by \({n}\). Hence, on the component indexed by \({n}\), \({\GmCan}\) is equal to \({\lambda_{n}(\mathbb{G}_m)\times B\mathbb{G}_m\subset \mathcal{I}_{B\mathbb{G}_m}}\), which is \({\mathcal{I}_{B\mathbb{G}_m}=\mathbb{G}_m\times B\mathbb{G}_m}\) if \({n\neq 0}\) and the trivial subgroup of \({\mathcal{I}_{B\mathbb{G}_m}}\) if \({n=0}\).
    In conclusion, we get that 
    \[\Grad^{\operatorname{rig}}=B\mathbb{G}_m\cup \coprod_{n\in \mathbb{Z}\setminus \{0\}}S\]
    i.e. for each connected component \({\mathfrak{Z}_n=B\mathbb{G}_m\subset \Grad}\) corresponding to \({\lambda_n:\mathbb{G}_m\to \mathbb{G}_m}\), \({\Zrig_n=S}\) if \({n\neq 0}\) and \({\Zrig_0=\mathfrak{Z}_0=B\mathbb{G}_m}\).
\end{example}

\subsection{Proof of the main theorem}

Here we use the action on the fibers of \({\gr:\Filt\to \Grad}\) together with \Cref{multip-act} to prove the remaining parts of \Cref{Theorem A}.

We begin with part 2:

\begin{theorem}\label{Part 2 of theorem A}
    Let \({\mathcal{X}}\) be an algebraic stack locally of finite presentation, with separated and quasi-compact diagonal and affine stabilizers over a quasi-separated scheme \({S}\). Let \({\mathfrak{Z}\subset \Grad}\) be a stack of gradings of \({\mathcal{X}}\) such that \({\gr:\FiltZ\to \mathfrak{Z}}\) has flat relative inertia of finite presentation over \({\FiltZ}\). Then the morphism
    \[\grRig:\FiltZRig\to \mathfrak{Z}\]
    is affine of finite presentation. In particular, if \({\gr:\FiltZ\to \mathfrak{Z}}\) is representable, it is affine.
\end{theorem}
\begin{proof}
    Let \({\rho:\FiltZ\to \FiltZRig}\) be the rigidification morphism.
    It is enough to show that \({\grRig:\FiltZRig\to \mathfrak{Z}}\) is affine after pullback along an atlas \({T\to \mathfrak{Z}}\) of \({\mathfrak{Z}}\).
    By \Cref{action on fibers 2} there is a contracting \({\mathbb{G}_m}\)-action on \({\FiltZRig\times_{\mathfrak{Z}}T}\), with fixed points given by \({\rho \circ \sigma_T(T)\cong T}\) by \Cref{fixed points in FiltRig} where \({\sigma_T:T\to \FiltZ\times_{\mathfrak{Z}}T}\) is the pullback of \({\sigma:\mathfrak{Z}\to \FiltZ}\) sending a grading to its associated split filtration. By part 1 of \Cref{multip-act}, \({\rho \circ \sigma_T(T)\cong T}\) is a closed subspace.
    This implies that the following diagram commutes
    \begin{equation}\label{gr affine digram}
    \begin{tikzcd}
	{(\FiltZRig \times_{\mathfrak{Z}}T)^0} & {\FiltZRig \times_{\mathfrak{Z}}T} \\
	& T
	\arrow["\shortmid"{marking}, hook, from=1-1, to=1-2]
	\arrow["\plim"', shift right, curve={height=5pt}, from=1-2, to=1-1]
	\arrow["\grRig", from=1-2, to=2-2]
	\arrow["\cong", from=2-2, to=1-1]
\end{tikzcd}
\end{equation}
    where \({\plim:\FiltZRig\times_\mathfrak{Z} T\to (\FiltZRig\times_\mathfrak{Z} T)^0}\) acts by \({f\mapsto 0.f}\), defined by the canonical extension of the contracting \({\mathbb{G}_m}\)-action to an \({(\mathbb{A}^1,\cdot)}\)-monoid action.
    By part 2 and 5 of \Cref{multip-act}, 
    \[{\plim:\FiltZRig\times_\mathfrak{Z} T\to (\FiltZRig\times_\mathfrak{Z} T)^0}\] 
    is affine and of finite presentation. The commutativity of diagram \ref{gr affine digram} implies the same holds for 
    \[\grRig:\FiltZRig\times_\mathfrak{Z} T\to T\] finishing the proof. 
 \end{proof}

 Lastly, we turn to part 3 of the main theorem, which we prove with a method analogous to the one used in \Cref{Stacky BB} to show that the attracting locus of a \(\mathbb{G}_m\)-action with its fixed points removed admits a geometric quotient projective over the fixed points.

\begin{theorem}\label{Part 3 of theorem A}
    Let \({\mathcal{X}}\) be an algebraic stack locally of finite presentation, with separated and quasi-compact diagonal and affine stabilizers over a quasi-separated scheme \({S}\). Let \({\mathfrak{Z}\subset \Grad}\) be a stack of gradings of \({\mathcal{X}}\) such that \({\gr:\FiltZ\to \mathfrak{Z}}\) has flat relative inertia of finite presentation over \({\FiltZ}\). Then the inclusion
    \[\FiltZnonsplitRig\subset \FiltZRig\]
    is an open immersion and the morphism 
    \[\FiltZnonsplitRig \to \Zrig\]
    admits a tame relative moduli space, smooth locally projective over \({\Zrig}\).    
\end{theorem}
\begin{proof}
    Since \({\rho\circ\sigma:\mathfrak{Z}\to \FiltZRig}\) is a section of the affine morphism \[{{\grRig:\FiltZRig\to \mathfrak{Z}}}\] 
    it is a closed immersion, hence \({\FiltZnonsplitRig\subset \FiltZRig}\) is an open immersion.

    As the statement we want to show is local on \({\Zrig}\), we may assume that \({\Zrig}\) is connected, with a connected atlas \({f_0:T\to \Zrig}\). If \({\GmCan}\) is trivial over \({\mathfrak{Z}}\), then \({\mathfrak{Z}=\Zrig}\), and by \Cref{canonical action and canonical subgroup} the action \({\mu^{\operatorname{can}}:\mathbb{G}_m\times \Filt_T\to \Filt_T}\) is trivial. Hence, by \Cref{action comparison}, \({\rsInducedFnullRig:\mathbb{A}^1\times\Filt_T\to \Filt_T}\) is the trivial action. It follows that the immersion \({\rho\circ\sigma:\mathfrak{Z}\to \FiltZRig}\) is an isomorphism, so \({\FiltZnonsplitRig}\) is empty and there is nothing to prove.
    
    Assume that \({\GmCan}\) is non-trivial over \({\mathfrak{Z}}\), then by \Cref{gmcan properties}, \({\mathfrak{Z}\to \Zrig}\) is a \({\mathbb{G}_m}\)-gerbe. 
    By \Cref{Part 2 of theorem A}, \({\FiltZRig\to \mathfrak{Z}}\) is the relative spectrum of a finitely presented quasi-coherent \({\mathcal{O}_{\mathfrak{Z}}}\)-algebra \({\mathcal{A}}\).
    Let \({p:T\to \Zrig}\) be a smooth atlas. Consider the pullback diagram 
    \[\begin{tikzcd}
	{\relSpec_{\mathfrak{Z}_T}p_{\mathfrak{Z}}^*\mathcal{A}} & {\relSpec_{\mathfrak{Z}}\mathcal{A}} \\
	{\mathfrak{Z}_T} & {\mathfrak{Z}} \\
	T & {\mathfrak{Z}^{\operatorname{rig}}}
	\arrow["{\operatorname{gr}}", from=1-2, to=2-2]
	\arrow[from=2-2, to=3-2]
	\arrow["p"', from=3-1, to=3-2]
	\arrow[from=2-1, to=3-1]
	\arrow["{p_{\mathfrak{Z}}}"', from=2-1, to=2-2]
	\arrow[from=1-1, to=1-2]
	\arrow[from=1-1, to=2-1]
\end{tikzcd}.\]
    Since \({\mathfrak{Z}_T}\) is a \({\mathbb{G}_m}\)-gerbe over \({T}\), it becomes trivial after passing to an étale cover, which again will be an atlas of \({\mathfrak{Z}}\), we may therefore assume \({{\mathfrak{Z}_T= B\mathbb{G}_{m,T}}}\). 
    Pulling back the affine morphism \({{\relSpec_{\mathfrak{Z}_T}p_{\mathfrak{Z}}^*\mathcal{A}\to B\mathbb{G}_{m,T}}}\) along the atlas \({T\to B\mathbb{G}_{m,T}}\) shows that \({{\relSpec_{\mathfrak{Z}_T}p_{\mathfrak{Z}}^*\mathcal{A}\cong [\relSpec_{T} A/\mathbb{G}_{m,T}]}}\) over \({B\mathbb{G}_{m,T}}\), where \({A=\bigoplus_{d\in \mathbb{Z}} A_d}\) is an \({\mathcal{O}_T}\)-algebra with a \({\mathbb{G}_{m,T}}\)-action defining the grading.
    
    By \Cref{canonical action and canonical subgroup} and \Cref{action comparison}, the grading of \({A}\) is up to a positive multiple the grading coming from the contracting action \({\rsInducedFnullRig:\mathbb{A}^1\times\relSpec_{\mathfrak{Z}_T}p_{\mathfrak{Z}}^*\mathcal{A}\to \relSpec_{\mathfrak{Z}_T}p_{\mathfrak{Z}}^*\mathcal{A}}\), and thus \({A_d=0}\) for all \({d<0}\).
   
    Furthermore, by \Cref{fixed points in FiltRig}, the essential image of \(\rho\circ \sigma \), is identified with the fixed points of the action on \(\FiltZRig\) over \(\mathfrak{Z}_T\). That is, \({\rho\circ \sigma (T)=[\relSpec_{T} (A/A_{>0})/\mathbb{G}_{m,T}]}\) where \(A_{>0}\) is the ideal sheaf of elements of positive weight in \(A\). Hence we can make the identification 
    \[\FiltZnonsplitRig\times_{\Zrig}T\cong [\relSpec_T A/\mathbb{G}_{m,T}]\setminus [\relSpec_T (A/A_{>0})/\mathbb{G}_{m,T}],\]
    and
    \[\FiltZnonsplitRig\times_{\Zrig}T\to \text{Proj}_TA\]
    is a coarse moduli space by \cite[Proposition 1.6.1 (iii)]{WeightedBlowUps}. 
    
    As \({A}\) is an algebra of finite presentation over \({\mathcal{O}_T}\), \({\text{Proj}_TA}\) is projective over \(T\).
    Moreover, as \({\FiltZnonsplitRig\times_{\Zrig}T}\) is a quotient stack of a \({\mathbb{G}_m}\)-action with finite stabilizers, it is a tame stack by \cite[Theorem 3.2]{Tame-stacks}, so \({{\FiltZnonsplitRig\times_{\Zrig}T\to \text{Proj}_TA}}\) is a tame moduli space.
    
    Since the relative inertia of \({\FiltZnonsplitRig\times_{\Zrig}T\to T}\) is finite and formation of relative inertia commutes with base change, the relative inertia of \({\FiltZnonsplitRig\to \Zrig}\) is finite. 
    Hence by \cite[Theorem 3.1]{TwistedStableMapsToTameStacks}, \({\FiltZnonsplitRig\to \Zrig}\) admits a relative coarse moduli space, 
    \({\overbar{\FiltZnonsplitRig}}\), over \({\Zrig}\). Moreover, \({\FiltZnonsplitRig\to \Zrig}\) is tame since tameness of a morphism is fppf local on the target by \cite[Definition 3.3]{TwistedStableMapsToTameStacks} and we saw above that \({\FiltZnonsplitRig\times_{\Zrig}T\to T}\) is tame.   
    In other words, \({\FiltZnonsplitRig\to \Zrig}\) admits a relative tame moduli space which is smooth and smooth locally projective over \({\Zrig}\).
\end{proof}

\begin{corollary}\label{Theorem A when Zrig space}
Let \({\mathcal{X}}\) be an algebraic stack locally of finite presentation, with separated and quasi-compact diagonal and affine stabilizers over a quasi-separated scheme \({S}\). Let \({\mathfrak{Z}\subset \Grad}\) be a stack of gradings of \({\mathcal{X}}\) such that the stack of filtrations centered at \({\mathfrak{Z}}\), \({\FiltZ}\), has flat relative inertia. 
Suppose furthermore that \({\Zrig}\) is an algebraic space and that \({\mathfrak{Z}\to \Zrig}\) is isomorphic to \({\Zrig\times B\mathbb{G}_m\to \Zrig}\), then the tame moduli space of 
\[\FiltZnonsplitRig\]
is projective over \({\Zrig}\).
\end{corollary}
\begin{proof}
    In this setting, the proof of \Cref{Part 3 of theorem A} goes through as above but we do not need to pass to an atlas of \({\Zrig}\) to describe the tame moduli as the projective spectrum of a finitely generated algebra over \({\Zrig}\).
\end{proof}

\section{Non-reductive quotients}
\label{non-reductive-quotients}
Here we apply \Cref{Theorem A} to deduce \Cref{Uhat intro}, i.e our version of \cite[Theorem 0.2]{U-hat} and  \cite[Remark 7.12]{bérczi2020projective}, that works over any Noetherian base scheme as well as some other corollaries on quotients by graded unipotent groups.

Before we get into these corollaries, we recall 
an example of the formation of quotients under an action of the graded unipotent group from \Cref{graded unip group example}, due to Hamilton, Hoskins and Jackson.

\begin{example}\cite{hamilton2024affine}\label{example over field}
    Let \(\Uhat\) be the graded unipotent group of upper triangular matrices in \(\operatorname{SL}_{2,k}\) over a field \(k\) and let \(V=\operatorname{Mat}_{2\times 2}\times \mathbb{A}^1\) be the vector space scheme of pairs of \(2\times 2\) matrices and a scalar and let \(X=\mathbb{P}(V)\).
    Let \(\Uhat\) act on \(V\) by conjugation on \(\operatorname{Mat}_{2\times 2}\) and trivially on \(\mathbb{A}^1\), as this action commutes with scalar multiplication, it descends to an action on \(X\). 
    When we restrict the action of to the diagonal \(\mathbb{G}_m\subset \Uhat\), it is described in homogeneous coordinates by
    \[t.[\begin{pmatrix}
        a & b\\ c & d
    \end{pmatrix}, w]=[\begin{pmatrix}
        a & t^2b\\ t^{-2}c & d
    \end{pmatrix}, w]=[\begin{pmatrix}
        t^2a & t^4b\\ c & t^2d
    \end{pmatrix}, t^2w].\]
    so one of the fixed point components is \({Z_0=\{c\neq 0, a=b=d=w=0\}\cong \Spec k}\), which has associated attracting locus \(X_0^+=\{c\neq 0\}\).
    The action of the unipotent radical of \(\Uhat\), \(\mathbb{G}_a\), restricted to \(X_0^+\) is free.
    In this case one can compute the quotients explicitly, namely 
    \[X_0^+/\mathbb{G}_a\cong \mathbb{A}^3\]
    which of course is affine over \(Z_0=\Spec k\). One way to see this is to notice that each \(\mathbb{G}_a\) orbit of \(X_0^+\) has a unique representative of the form
    \[[\begin{pmatrix}
        a & b\\ 1 & 0
    \end{pmatrix}, w].\]
    The action of \(\mathbb{G}_m\) on \(X_0^+\) descends to an action on \(\mathbb{A}^3\), described by \({t.(a,b,w)=(t^2a,t^4b,t^2w)}\), and \(Z_0\) is identified with the fixed point \((0,0,0)\) under the quotient morphism \(X^+_0\to \mathbb{A}^3\), so letting \(X_0^{+,\circ}:=X_0^{+}\setminus Z\), we have \(X_0^{+,\circ}/U\cong \mathbb{A}^3\setminus\{0\}\). We conclude that
    \[X_0^{+,\circ}/\Uhat\cong \mathbb{P}(1,1,2),\]
    that is, the \(\Uhat\)-quotient is isomorphic to the weighted projective plane over \(k\) with weights \((1,1,2)\).
\end{example}

\subsection{General graded unipotent actions}
The main point of this section is to use \Cref{Theorem A} to deduce a theorem about general algebraic spaces with an action of a graded unipotent group, which we will use in the next section to prove \Cref{Uhat intro}.

Let \(X\) be an algebraic space with an action of a graded unipotent group \(\Uhat=U\rtimes \mathbb{G}_m\). Let \(X^{+}\) and \(X^{0}\) be the spaces of attractors and fixed points for the \(\mathbb{G}_m\)-action on \(X\) induced by \(\mathbb{G}_m\subset \Uhat\) and fix a subspace \(Z\subset X^{0}\). Let \(\AttractZ\), the attractor centered at \(Z\), be defined via the cartesian diagram
\[\begin{tikzcd}
	\AttractZ & {X^+} \\
	Z & {X^{0}}
	\arrow[from=1-1, to=1-2]
	\arrow[from=1-1, to=2-1]
	\arrow["\plim", from=1-2, to=2-2]
	\arrow[from=2-1, to=2-2]
\end{tikzcd}\]
As noted in \Cref{filt of quotients in text}, there is a natural action of \(\Uhat\) on \(\AttractZ\).\\
Using this action, let \(\AttractZnonsplit:=\AttractZ\setminus \Uhat.Z\) be the complement of the \(\Uhat\)-sweep of \(Z\subset \AttractZnonsplit\), which is stable under the action of \(\Uhat\) on \(\AttractZ\).  

\begin{theorem}\label{Theorem B in text}
        Let \(X\) be a quasi-separated algebraic space locally of finite type over a Noetherian scheme \(S\) with an action of a graded unipotent group \(\Uhat\).  Let \(Z\subset X^{0}\) be a subspace of the \(\mathbb{G}_m\)-fixed points of \(X\). If the stabilizer scheme of \(U\) acting on \(\AttractZ\) is flat over \(\AttractZ\), then
    \begin{enumerate}
        \item the algebraic space \(\AttractZ\) admits a geometric \(U\)-quotient, which is affine and of finite type over \(Z\).
        \item the algebraic space \(\AttractZnonsplit\) admits a geometric \(\Uhat\)-quotient, projective over \(Z\).
    \end{enumerate}
\end{theorem}

\begin{proof}\label{Non-reductive-quotients}
    By \Cref{Filt of quotients} we have a cartesian diagram 
    \[\begin{tikzcd}
	{[\AttractZ/\Uhat]} & {\operatorname{Filt}([X/\Uhat])} \\
	{Z\times B\mathbb{G}_m} & {\operatorname{Grad}([X/\Uhat])}
	\arrow[from=1-1, to=1-2]
	\arrow["{P_{\lim}}"', from=1-1, to=2-1]
	\arrow["\gr", from=1-2, to=2-2]
	\arrow[from=2-1, to=2-2]
\end{tikzcd}\]
    with the horizontal arrows being open immersions. 
    
    We first explain how the conclusion of the theorem holds if the relative inertia of\\ \(P_{\lim}:[\AttractZ/\Uhat]\to Z\times B\mathbb{G}_m\), \(\mathcal{I}_{P_{\lim}}\) is flat, and then show that our assumptions on the \(U\)-stabilizers imply that \(\mathcal{I}_{P_{\lim}}\) is flat.\footnote{The finite presentation assumption of the relative inertia is automatic as we are working with quotient stacks under actions of finitely presented group schemes.} 
    
    Letting \(\mathcal{X}=[X/\Uhat]\) and \(\mathfrak{Z}:=Z\times B\mathbb{G}_m\subset \Grad\), \Cref{Theorem A} now implies that, if \(\mathcal{I}_{P_{\lim}}\) is flat over \([\AttractZ/\Uhat]\), then \([\AttractZ/\Uhat]\to Z\times B\mathbb{G}_m\) factors as a gerbe and an affine morphism. 
    In particular, the morphism \([\AttractZ/\Uhat]\to Z\times B\mathbb{G}_m\) admits a relative coarse moduli space, affine over \(Z\times B\mathbb{G}_m\). 
    Pulling back \(P_{\lim}:[\AttractZ/\Uhat]\to Z\times B\mathbb{G}_m\) along the atlas \(Z\to Z\times B\mathbb{G}_m\), this implies that \([\AttractZ/U]\to Z\) admits a relative coarse moduli space, affine over \(Z\). In other words, \(\AttractZ\) admits a geometric \(U\)-quotient affine over \(Z\), giving part 1.
    
    With notation as in the previous paragraph, \(\Zrig=Z\), and \(\FiltZnonsplit=[\AttractZnonsplit/\Uhat]\) by \Cref{Non-split filt in quotients}. Hence if \(\mathcal{I}_{P_{\lim}}\) is flat, then \([\AttractZnonsplit/\Uhat]\to Z\) admits a relative coarse moduli space, projective over \(Z\) by \Cref{Theorem A when Zrig space}. In other words, if \(\mathcal{I}_{P_{\lim}}\) is flat over \([\AttractZ/\Uhat]\), then \(\AttractZnonsplit\) admits a geometric \(\Uhat\)-quotient, projective over \(Z\), giving part 2 of the theorem.

    To finish the proof we need to show that \(\text{Stab}_{U,\AttractZ}\) being flat over \(\AttractZ\) implies that \(\mathcal{I}_{P_{\lim}}\) is flat over \([\AttractZ/\Uhat]\).

    The relative inertia \(\mathcal{I}_{P_{\lim}}\) is the kernel of the natural morphism 
    \[\mathcal{I}_{[\AttractZ/\Uhat]}\to P_{\lim}^*\mathcal{I}_{Z\times B\mathbb{G}_m}=\mathbb{G}_{m,[\AttractZ/\Uhat]}\]
    of group stacks over \([\AttractZ/\Uhat]\). This morphism pulls back along the atlas \(\AttractZ\to [\AttractZ/\Uhat]\) to a map 
    \[\text{Stab}_{\Uhat,\AttractZ}\to \mathbb{G}_{m,\AttractZ}\]
    defined by \((x,ut)\mapsto (x,t)\) for \(ut\in \Uhat\) stabilizing \(x\), using that \(\text{Stab}_{\Uhat,\AttractZ}\) is a subgroup scheme of \(\AttractZ\times \Uhat\) and the fact that \(\Uhat=U\rtimes \mathbb{G}_m\), \Cref{semi-dir decomp}. 
    In particular, this map takes \((x,ut)\) to \((x,1)\in \mathbb{G}_{m,\AttractZ}\) if and only if \(t=1\), so the kernel is \(\text{Stab}_{U,\AttractZ}\). In particular, the kernel is flat by assumption and hence so is \(\mathcal{I}_{P_{\lim}}\) by descent.
    \end{proof}

\begin{corollary}\label{Contracting Theorem B}
    Let \(X\) be a quasi-separated algebraic space locally of finite type over a Noetherian scheme \(S\) with an action of a graded unipotent group \(\Uhat\). If the action of the grading \(\mathbb{G}_m\subset \Uhat\) is contracting and the stabilizer group scheme for the \(U\)-action on \(X\) is flat, then 
     \begin{enumerate}
        \item  \(X\) admits a geometric \(U\)-quotient, which is affine and of finite type over the \(\mathbb{G}_m\)-fixed points \(X^{0}\).
        \item  \(X\setminus U.X^{0}\) admits a geometric \(\Uhat\)-quotient, projective over \(X^{0}\).
    \end{enumerate}
\end{corollary}
\begin{proof}
    By \Cref{multip-act} (5), \(X\cong X^+\), hence the statement follows directly from \Cref{Theorem B in text}.
\end{proof}

\subsection{Linear actions on projective schemes}
Here we use \Cref{Contracting Theorem B} to deduce our version of the \(\Uhat\)-theorem, \Cref{Uhat intro}.

Let \(X\) be a reduced projective scheme with irreducible geometric fibers over an affine Noetherian scheme \(S\) with an action of a graded unipotent group scheme over \(S\), \(\hat{U}=U\rtimes \mathbb{G}_m\). 

Suppose this action admits a very ample linearization \(L\)\footnote{That is a lift of the action on \(X\) to the line bundle, inducing an equivariant closed immersion \({X\to  \mathbb{P}(H^{0}(X,L)^{\lor})}\), where \(V^{\lor}\) denotes the dual of a module \(V\).}, such that the action of \(\mathbb{G}_m\)  on \(H^{0}(X,L)^{\lor}\) via the inclusion \(\mathbb{G}_m\subset \Uhat\) has weights
\[w_{\text{min}}=w_0<w_1<\ldots<w_h=w_{\text{max}}.\]

Let \(Z_{\text{min}}\subset X^0\) be the subspace of fixed points where the action of \(\mathbb{G}_m\subset \Uhat\) on the stalks of \(L^{\lor}\) acts by weight \(w_{\text{min}}\), and let \(X_{\text{min}}^+:=X^+_{Z_{\text{min}}}\) be the attractor centered at \(Z_{\text{min}}\). 
\begin{definition}\label{adapted and well-adapted linearizations}
    We say that the linearization \(L\) of the graded unipotent action is
    \begin{itemize}
        \item \textit{adapted} if \(w_{\text{min}}<0<w_1\),
        \item \textit{well-adapted} if it is adapted, and for all \(x\in X\setminus X^{\mathbb{G}_m}\), \(\operatorname{Stab}_{\mathbb{G}_m}\) acts trivially on the stalk \(L_x\).
    \end{itemize}
\end{definition}
A computation in the coordinates given by the linearization shows that \({\evone:X_{\text{min}}^+\to X}\) is an open immersion onto the open Białynicki-Birula stratum for the \(\mathbb{G}_m\)-action on \(X\).\\
Let \(X^{ss,U}:=X_{\text{min}}\setminus \Uhat.Z_{\text{min}}\), if the action is adapted, then \(X^{ss,U}\) is the Hilbert-Mumford stable locus of the linearization. 
By \Cref{Filt of quotients} and \Cref{Non-split filt in quotients}, if follows that
\[[X^{ss,U}/\Uhat]=\text{Filt}([X/\Uhat])^{\circ}_{[Z_{\text{min}}/\mathbb{G}_m]}\]
is the stack of non-split filtrations of \([X/\Uhat]\) centered at \([Z_{\text{min}}/\mathbb{G}_m]\).

\begin{theorem}\label{Classical Uhat}
Let \(X\) be a reduced projective scheme with irreducible geometric fibers over a Noetherian affine scheme \(S\), with an action of a graded unipotent group \(\hat{U}\), linearized by a very ample line bundle \(L\).
Suppose the action of \(U\) on \(X_{\operatorname{min}}\) has smooth stabilizers with constant dimension, then 
\begin{enumerate}
    \item the action of \(U\) on \(X_{\operatorname{min}}\) admits a geometric quotient affine over \(Z_{\operatorname{min}}\)
    \item the action of \(\hat{U}\) on \(X^{ss,U}\) admits a geometric quotient \(X^{ss,U}/\hat{U}\), projective over \(S\).
\end{enumerate}
If the action of \(U\) is free, the above quotients are tame, i.e. pushforward of equivariant quasi-coherent sheaves along the quotient morphism is exact. If furthermore the linearization is well adapted, then \(L\) descends to a line bundle on \(X^{ss,U}/\hat{U}\).

\end{theorem}
\begin{remark}
    The version of the theorem stated in the introduction, \Cref{Uhat intro}, follows from the observation that \(X_{\operatorname{min}}\) is the open Białynicki-Birula stratum of \(X\) for the grading \(\mathbb{G}_m\)-action. 
\end{remark}
\begin{proof}
    By assumption \(X\) is reduced, hence the open subscheme \(X_{\text{min}}\) is as well.
    We have also assumed that \(U\) acts with smooth stabilizers of constant dimension on \(X_{\text{min}}\), hence the stabilizer group scheme \(\text{Stab}_U\to X_{\text{min}}\) is smooth  by \cite[Tome VIb, Corollaire 4.4]{SGA3Vol2} which states that a group scheme of finite type over a reduced scheme with smooth fibers of constant dimension is smooth. In particular, \(\text{Stab}_U\to X_{\text{min}}\) is flat. 
    Thus the assumptions of \Cref{Contracting Theorem B} are satisfied so we get a geometric quotient for the \(U\)-action on \(X_{\text{min}}\), affine over \(Z_{\text{min}}\), and one of the \(\Uhat\)-action on \(X^{ss}\), projective over \(Z_{\text{min}}\). Since \(Z_{\text{min}}\) is a closed subscheme of \(X\), a projective scheme, \(Z_{\text{min}}\) is projective over \(S\) and hence so is \(X^{ss}/\Uhat\).

    If the \(U\)-action is free, then \(\FiltQZ=\FiltQZRig\) and it follows that \({X^{ss,U}\to X^{ss,U}/\Uhat}\) is a tame quotient.  
    Let \({\pi:[X^{ss,U}/\Uhat]\to X^{ss,U}/\Uhat}\) be the tame moduli morphism.
    The line bundle \(L\) descends to a line bundle \(\mathcal{L}'\) on \([X^{ss,U}/\Uhat]\) as it linearizes the action of \(\Uhat\) on \(X\). Since the linearization is assumed to be well adapted, it follows that the inertial action on stalks of \(\mathcal{L}'\) for points in \([X^{ss,U}/\Uhat]\) is trivial. Hence, by \cite[Theorem 10.3]{alperGMS}, \(\mathcal{L}'\) descends to a line bundle \(\mathcal{L}=\pi_*\mathcal{L}'\) on \(X^{ss,U}/\Uhat\).

\end{proof}

\begin{remark}[Comparison with \cite{U-hat}]\label{comparison with classical uhat}
    Given any ample linearization of an action of \(\Uhat\) on \(X\) with at least two distinct weights, one can modify it to one that is adapted by twisting by a character of \(\Uhat\). By taking tensor powers of \(L\) and twisting by a character, one can obtain a well adapted linearization. The definition of adapted and well adapted linearizations in \cite{U-hat} differs from our characterization in that they are defined in terms of these twisting characters. 
    
    In \cite[Theorem 0.2]{U-hat} it is shown, in the case \(S=\Spec \mathbb{C}\), that if the linearization is well adapted, then the line bundle descended to \(X^{ss,U}/\Uhat\) is ample and furthermore, that
    \[X^{ss,U}/\Uhat\cong \operatorname{Proj}_{\mathbb{C}}\left(\bigoplus_{n\geq 0}H^0(X,L^n)\right)^{\Uhat}\]
    in particular showing that the invariant algebra \(\left(\bigoplus_{n\geq 0}H^0(X,L^n)\right)^{\Uhat}\) is finitely generated. Currently I don't see how to show this using our methods, but I do believe it holds in the generality of \Cref{Classical Uhat}.
\end{remark}

\begin{appendix}
\section{Filtrations of quotient stacks}
In this appendix we describe stacks of gradings and filtrations of quotient stacks, as quotient stacks. The main result is a version of  \cite[Theorem 1.4.8]{halpernleistner2022structure}. 
The description of \cite[Theorem 1.4.8]{halpernleistner2022structure} works for actions of algebraic groups over a field with split maximal tori and describes the whole stacks of filtrations and gradings in terms of cocharacters of the split maximal torus, which up to conjugacy are defined over the base field. If there is no split maximal torus, we can not give the same global description but we can still describe some filtrations with a particular center as long as \(G\) admits some cocharacters defined over the base.

The goal of this appendix is to show the following.

\begin{theorem}\label{Filt of quotients}
    Let \(G\) be a smooth connected affine group scheme over a Noetherian scheme \(S\) with an action on a quasi-separated algebraic space \(X\) locally of finite type over \(S\) and let \(\lambda:\mathbb{G}_{m}\to G\)  be a cocharacter. The diagram 
    \[\begin{tikzcd}
	{[X^{+,\lambda}/P_{\lambda}]} & {\operatorname{Filt}([X/G])} \\
	{[X^{0,\lambda}/L_{\lambda}]} & {\operatorname{Grad}([X/G])}
	\arrow["F_{\lambda}", from=1-1, to=1-2]
	\arrow["P_{\lim}"',from=1-1, to=2-1]
	\arrow["G_{\lambda}",from=2-1, to=2-2]
	\arrow["\gr",from=1-2, to=2-2]
\end{tikzcd}\]
with maps defined as in \Cref{filt of quotients in text}, commutes and the horizontal arrows are open immersions.
\end{theorem}

It is shown in \cite[Theorem 1.4.7]{halpernleistner2022structure} that this follows if one proves the theorem in the special case \([X/G]=BG\), which is what we will do below. 

To describe the stack of gradings of a classifying stack \(BG\), we will make use of the sheaf of cocharacters of \(G\),
\[\text{Hom}_{gp}(\mathbb{G}_m,G).\]
This sheaf is represented by a smooth separated scheme over \(S\) if \(G\) is smooth and affine over \(S\) by \cite[Expose XI, Corollarie 4.2]{SGA3Vol2} with a \(G\)-action by conjugation. 

The following lemma is a corollary to \cite[Lemma 1.1.5]{halpernleistner2022structure}, which states that for any stack \(\mathcal{X}\) and any \(\xi:S\to\mathcal{X}\), the diagram
\[\begin{tikzcd}
	{\text{Hom}_{gp}(\mathbb{G}_m,\text{Aut}_{\mathcal{X}}(\xi))} & S \\
	{\operatorname{Grad}(\mathcal{X})} & \mathcal{X}
	\arrow[from=1-1, to=2-1]
	\arrow[from=1-1, to=1-2]
	\arrow["\xi", from=1-2, to=2-2]
	\arrow["u"', from=2-1, to=2-2]
\end{tikzcd}\]
is cartesian, where \(u:\operatorname{Grad}(\mathcal{X})\) is the morphism mapping a \(T\)-point \(g\in \Grad\) defined by \(B\mathbb{G}_{m,T}\to \mathcal{X}\) to the composition \(T\to B\mathbb{G}_{m,T}\to \mathcal{X}\) where the first morphism is the atlas of \(B\mathbb{G}_{m,T}\). The morphism \({\text{Hom}_{gp}(\mathbb{G}_m,\text{Aut}_{\mathcal{X}}(\xi))\to \Grad}\) is defined by mapping a cocharacter \(\lambda:\mathbb{G}_{m,T}\to \text{Aut}_{\mathcal{X}}(\xi))\) defined over a \(S\)-scheme \(T\) to the morphism \(f_{\xi,\lambda}:B\mathbb{G}_m\to \mathcal{X}\) defined by \(\xi:T\to \mathcal{X}\) and \(\lambda\). Two such morphisms \(f_{\xi,\lambda_1}\) and \(f_{\xi,\lambda_2}\) are isomorphic in \(\Grad\) if and only if \(\lambda_1\) and \(\lambda_2\) are conjugate via some \(g\in \text{Aut}_{\mathcal{X}}(\xi)(T)\).

\begin{lemma}\label{Grad of classifying}
    Let \(G\) be a connected smooth affine group scheme over \(S\), then 
    \[\operatorname{Grad}(BG)\cong [\text{Hom}_{\text{gp}}(\mathbb{G}_m,G)/G]\]
    where \(G\) acts by conjugation.
    In particular, if \(\lambda:\mathbb{G}_{m,S}\to G\) is a cocharacter defined over \(S\), then 
    \(BL_{\lambda}\) is an open substack of \(\operatorname{Grad}(BG)\) and if a torsor over \(B\mathbb{G}_{m,T}\) defines a point in \(BL_{\lambda}(T)\subset \operatorname{Grad}(BG)(T)\), it admits a canonical reduction of structure group to \(L_{\lambda}\) after passing to an étale cover of \(T\).
\end{lemma}
\begin{proof}
By \cite[Lemma 1.1.5]{halpernleistner2022structure}, the following diagram is cartesian    
\[\begin{tikzcd}
	{\text{Hom}_{gp}(\mathbb{G}_m,G)} & S \\
	{\operatorname{Grad}(BG)} & BG
	\arrow[from=1-1, to=2-1]
	\arrow[from=1-1, to=1-2]
	\arrow["p", from=1-2, to=2-2]
	\arrow["u"', from=2-1, to=2-2]
\end{tikzcd}\]
where \(u:\operatorname{Grad}(BG)\to BG\) is defined by mapping \(f:B\mathbb{G}_{m,T}\to BG\) to the precomposition by the atlas \(T\to B\mathbb{G}_{m,T}\) and
\(p:S\to BG\) is the standard atlas of \(BG\). 
In particular, \(\text{Hom}_{gp}(\mathbb{G}_m,G)\) is a \(G\)-torsor over \(\operatorname{Grad}(BG)\), and the action of \(G\) is by conjugation of cocharacters. The fact that the action is given by conjugation follows from the explicit description of the morphism \(\text{Hom}_{gp}(\mathbb{G}_m,G)\to \operatorname{Grad}(BG)\) given in \cite[Lemma 1.1.5]{halpernleistner2022structure} and briefly above, which proves the first part of the lemma.

Let \(\lambda:\mathbb{G}_{m,S}\to G\) be a cocharacter and \(L_{\lambda}\subset G\) its centralizer. By \cite[Expose XI, Corollaire 5.3]{SGA3Vol2}, \(G/L_{\lambda}\subset \text{Hom}_{gp}(\mathbb{G}_m,G)\) is a \(G\)-equivariant open immersion parameterizing the conjugacy class of \(\lambda\). This induces an open immersion 

\[BL_{\lambda}=[(G/L_{\lambda})/G]\subset [\text{Hom}_{gp}(\mathbb{G}_m,G)/G]=\operatorname{Grad}(BG).\]

We now show that any \(G\)-torsor \(\mathcal{P}\) over \(B\mathbb{G}_{m,T}\) defining a point \(f_0\in \Grad(T)\) factoring through \(BL_{\lambda}\) admits a canonical reduction of structure group to \(L_{\lambda}\) after pulling back to an étale cover. Note that the \(G\)-torsor on \(T\) defined by pulling back \(\text{Hom}_{gp}(\mathbb{G}_m,G)\to \operatorname{Grad}(BG)\) admits a section étale locally as \(G\) is smooth. Hence there is an étale covering \(T'\to T\) such that \(f_{0}|_{T'}:\Grad\) admits a lift to \(\text{Hom}_{gp}(\mathbb{G}_m,G)\), in particular \(f_{0}|_{T'}\) is isomorphic to the point of \(\mathcal{P}_{\lambda}\in\Grad(T')\) defined by \(\lambda:\mathbb{G}_{m,T'}\to G\). In particular, there is a \(T'\)-section of \(\mathcal{P}_{\lambda}\) which induces an isomorphism \(G\cong \text{Aut}(\mathcal{P}_{\lambda})\) and a section \(B\mathbb{G}_m\to \mathcal{P}_{\lambda}\), and in particular a section of \(\mathcal{P}_{\lambda}/L_{\lambda}\), that is a reduction of structure group to \(L_{\lambda}\).
\end{proof}

Let us now describe the filtrations of \(BG\) centered at \(BL_{\lambda}\), that is \(\FiltZ\) for \(\mathcal{X}=BG\) and \(\mathfrak{Z}=BL_{\lambda}\).

\begin{lemma}\label{Filt of classifying}
    Let \(G\) be a connected smooth affine group scheme over \(S\).
    For any cocharacter \(\lambda:\mathbb{G}_{m,S}\to G\), the diagram 
    \[\begin{tikzcd}
	{BP_{\lambda}} & {\operatorname{Filt}(BG)} \\
	{BL_{\lambda}} & {\operatorname{Grad}(BG)}
	\arrow["F_{\lambda}",from=1-1, to=1-2]
	\arrow["P_{\lim}"',from=1-1, to=2-1]
	\arrow["G_{\lambda}",from=2-1, to=2-2]
	\arrow["\gr",from=1-2, to=2-2]
\end{tikzcd}\]
is cartesian.
\end{lemma}
The proof goes along the same lines as the proof of  \cite[Lemma 1.7]{HM-Stacks}, describing \(G\)-torsors on \(\Theta\) for \(G\) reductive over an algebraically closed field. The idea is to show that given a \(T\)-point \(f_0\in \operatorname{Grad}(BG)\), all points in the fiber of \(\gr:\operatorname{Filt}(BG)\to \operatorname{Grad}(BG)\) admit a canonical reduction of structure group to \(P_{\lambda}\) when viewed as a \(G\)-torsor on \(\Theta\). This is shown using deformation theory and the fact that the embedding \(B\mathbb{G}_m\to \Theta\) is coherently complete, \cite[Proposition 5.1]{alper2023etale}.

\begin{proof}
    Let \(\mathfrak{Z}=BL_{\lambda}\subset \operatorname{Grad}(BG)\).
    We first show that \(F_{\lambda}:BP_{\lambda}\to \operatorname{Filt}(BG)\) factors through \(\operatorname{Filt}(BG)_{\mathfrak{Z}}\).
    Let \(f_{\lambda}:\Theta\to BG\) be the morphism defined by the structure morphism \(\mathbb{A}^1\to S\) and the cocharacter \(\lambda:\mathbb{G}_m\to G\), in other words \(f_{\lambda}\) is \(F_{\lambda}\) applied to the atlas \(S\to BP_{\lambda}\). 
    The restriction of \(f_{\lambda}:\Theta\to BG\) to \(B\mathbb{G}_m\) is the composition \(S\to \text{Hom}_{gp}(\mathbb{G}_m,G)\to \operatorname{Grad}(BG)\) where the first morphism is defined by \(\lambda\) and the second is the atlas of \(\operatorname{Grad}(BG)\) defined in \Cref{Grad of classifying}, i.e \(\gr(f_{\lambda})\in BL_{\lambda}\). Since all points of \(BP_{\lambda}\) étale locally factors through the atlas \(S\to BP_{\lambda}\), we get the factorization \(F_{\lambda}:BP_{\lambda}\to \operatorname{Filt}(BG)_{\mathfrak{Z}}\subset \operatorname{Filt}(BG)\). 

    To see that \(F_{\lambda}:BP_{\lambda}\to \operatorname{Filt}(BG)_{\mathfrak{Z}}\) is surjective, it is enough to show that all \(T\) points of \(\operatorname{Filt}(BG)_{\mathfrak{Z}}\) are étale locally isomorphic to \(F_{\lambda}(\mathcal{F})\) for some \(P_{\lambda}\) torsor \(\mathcal{F}\).
    Let \(\mathcal{P}\) be a \(G\)-torsor over \(\Theta_T\)  defining a point in \(\operatorname{Filt}(BG)_{BL_{\lambda}}\).
    As we may work étale locally on \(T\), we may assume that \(\mathcal{P}_{0}=\mathcal{P}|_{B\mathbb{G}_m}\) admits a canonical reduction of structure group to \(L_{\lambda}\), \(B\mathbb{G}_m \to \mathcal{P}_{0}/L_{\lambda}\) by \Cref{Grad of classifying}.
    In particular, there is canonical reduction \(s_0:B\mathbb{G}_m\to \mathcal{P}/P_{\lambda}\), as there is a surjection \(\mathcal{P}_{0}/L_{\lambda}\to \mathcal{P}_{0}/P_{\lambda}\).
    
    We now use deformation theory to show that this reduction lifts uniquely to a reduction of structure group of \(\mathcal{P}|_{\Theta^{[n]}}\), 
    where \(\Theta^{[n]}\) denotes the \(n\)-th infinitesimal thickening of \({B\mathbb{G}_{m,T}\subset \Theta_T}\), i.e. it is defined by the ideal sheaf \((x^{n+1})\) on the atlas. 
    
    Since \(s_0\) is a section of a representable morphism it is representable, and thus we can apply \cite[Theorem 1.5]{DeformationOfRepresentable} to see that the obstruction to lifting a section  \(s_n:\Theta^{[n]}\to \mathcal{P}/P_{\lambda}\) to \(s_{n+1}:\Theta^{[n+1]}\to \mathcal{P}/P_{\lambda}\) lies in 
    \[\text{Ext}^1(s_n^*L_{(\mathcal{P}/P_{\lambda})/\Theta_T},(x^{n+1})/(x^{n+2})).\] 
    Similarly if the obstruction vanishes, the isomorphism classes of deformations forms a torsor under 
     \[\text{Ext}^0(s_n^*L_{(\mathcal{P}/P_{\lambda})/\Theta},(x^{n+1})/(x^{n+2})).\] 
    However, the cotangent complex of \(\mathcal{P}\) over 
    \(\Theta_T\) is quasi-isomorphic to the complex with 
    \(\text{Lie}(G)^{\lor}\) in degree 0 with \(\mathbb{G}_m\) acting via \(\lambda\) and the dual of the adjoint representation, and the cotangent complex of
    \(\mathcal{P}/P_{\lambda}\) is quasi-isomorphic to the complex with 
    \((\text{Lie}(G)/\text{Lie}(P_{\lambda}))^{\lor}\) in degree 0, and thus \(\mathbb{G}_{m,T}\) acts on it with strictly positive weights\footnote{This is due to the fact that \(P_{\lambda}\subset G\) is the attractor for the conjugation action of \(\mathbb{G}_m\) on \(G\) through \(\lambda\).}. Hence
    \[\mathbb{R}\text{Hom}(s_n^*L_{(\mathcal{P}/P_{\lambda})/\Theta_{T}},(x^{n+1})/(x^{n+2}))\cong 0\]
    in the derived category of quasi-coherent sheaves on \(\Theta_{T}\), since the weight of \((x^{n+1})/(x^{n+2})\) is strictly negative.\footnote{This follows from the ideal \((x^{n+1})\) being a line bundle on \(\Theta\) with a section for all \(n\geq0\). See \cite[Section 1.A.a]{HM-Stacks}.}. Thus the \(\operatorname{Ext}\)-groups above vanish and there is a canonical sequence of infinitesimal lifts \(s_n\) of \(s_0\).
    The pair \((\Theta,B\mathbb{G}_m)\) is coherently complete by \cite[Theorem 1.6]{alper2023etale} and the fact that \(\Theta\to \operatorname{pt}\) is a good moduli space morphism. Hence the system of infinitesimal lifts induces a lift of the section \(s_0:B\mathbb{G}_m\to \mathcal{P}/P_{\lambda}\) to a a global one \(s:\Theta\to \mathcal{P}/P_{\lambda}\) . So there is a canonical reduction of structure group of \(\mathcal{P}\) to \(P_{\lambda}\).
    In particular, \({F_{\lambda}:BP_{\lambda}\to \operatorname{Filt}(BG)_{BL_{\lambda}}}\) is surjective. 
    
    Since the reduction of structure group is canonical, it follows that \(F_{\lambda}:BP_{\lambda}\to \FiltZ\) is also a monomorphism, finishing the proof.
\end{proof}

This finishes the proof of \Cref{Filt of quotients} in the case \([X/G]=BG\), and thus in general by the argument presented in \cite[\textit{The general case} in the proof of Theorem 1.4.7]{halpernleistner2022structure}.
\section{Universal property of rigidification}
Here we collect some properties of the rigidification of a stack with respect to a subgroup stack of its inertia stack.

Throughout this section \(\mathcal{X}\) is an algebraic stack locally of finite presentation over an algebraic space \(S\) and \(H\subset \mathcal{I}_{\mathcal{X}}\) a closed subgroup stack of the inertia of \(\mathcal{X}\), that is \(H\) is preserved under composition and inversion in \(\mathcal{I}_{\mathcal{X}}\) and the unit section \(\mathcal{X}\to \mathcal{I}_{\mathcal{X}}\) factors through \(H\). 
In particular, for each \(x:T\to \mathcal{X}\), we get a flat finitely presented subgroup \(H_x\subset \text{Aut}_T(x)\), which is in fact a normal subgroup as shown in \cite{Tame-stacks}. It is shown in \cite{Tame-stacks} that there is a morphism of locally finitely presented algebraic stacks \(\rho:\mathcal{X}\to\mathcal{X}\!\!\!\fatslash H\), unique up to equivalence, such that
\begin{enumerate}
    \item \(\rho:\mathcal{X}\to \mathcal{X}\!\!\!\fatslash H\) is a smooth \textit{fppf}-locally trivial gerbe and
    \item for each \(x:T\to\mathcal{X}\), the homomorphism of group schemes
    \[\text{Aut}_T(x)\to \text{Aut}_T(\rho\circ x)\]
    is surjective with kernel \(H_x\).
\end{enumerate}

\begin{lemma}[Formation of rigidification commutes with base change]\label{rigi commutes with base change}
    Let \(S'\to S\) be a morphism of algebraic spaces. The rigidification of \(\mathcal{X}_{S'}\) with respect to \(H_S'\) is isomorphic to \((\mathcal{X}\!\!\!\fatslash H)_S'\).
\end{lemma}
\begin{proof}
    By the characterizing properties of the rigidification, it is enough to show that 
    \begin{enumerate}
        \item \(\rho_{S'}:\mathcal{X}_{S'}\to(\mathcal{X}\!\!\!\fatslash H)_S'\) is a smooth fppf trivializable gerbe, and
        \item for all \((x,s):T\to \mathcal{X}_{S'}\), the morphism 
        \[\text{Aut}_T(x,s)\to \text{Aut}_T(\rho(x),s)\]
        is surjective with kernel \(H_{(x,s)}\). 
    \end{enumerate}
    Since being smooth and being an fppf trivializable gerbe both are properties stable under base change, it is enough to check the second point. This follows from \(\text{Aut}_T(x,s)\to  \text{Aut}_T(\rho(x),s)\) being the pullback of \(\text{Aut}_T(x)\to  \text{Aut}_T(\rho(x))\), which has the property by assumption.
\end{proof}

\begin{lemma}[Universal property]\label{universal property of rigi}
    Let \(f:\mathcal{X}\to \mathcal{Y}\) be a morphism of algebraic stacks and suppose \(H\subset \mathcal{I}_{f}\) is a closed subgroup stack of the relative inertia of \(f\), flat and of finite presentation over \(\mathcal{X}\), then there is a unique morphism 
    \[f^{rig}:\mathcal{X}\!\!\!\fatslash H\to \mathcal{Y}\] 
    such that 
    \[f=f^{rig}\circ \rho.\]
\end{lemma}
\begin{proof}
    We prove this by constructing an atlas of \(\mathcal{X}\!\!\!\fatslash H\) adapted to the morphism \(f:\mathcal{X}\to \mathcal{Y}\).
    Let \(Y\to \mathcal{Y}\) be an atlas and consider the diagram 
\[\begin{tikzcd}
	X & {\mathcal{X\times_{Y}}Y} & Y \\
	& {\mathcal{X}} & {\mathcal{Y}}
	\arrow[from=1-2, to=1-3]
	\arrow[from=1-3, to=2-3]
	\arrow[from=1-2, to=2-2]
	\arrow["f"', from=2-2, to=2-3]
	\arrow["p",from=1-1, to=2-2]
	\arrow[from=1-1, to=1-2]
\end{tikzcd}\]
    where \(Y\to \mathcal{Y}\) and \(X\to \mathcal{X\times_{Y}}Y\) are atlases, so \(p:X\to \mathcal{X}\) is an atlas as well. Then the morphism
    \[f:\mathcal{X}\to \mathcal{Y}\] 
    induces a morphism of groupoids
\[\begin{tikzcd}
	{X\times_{\mathcal{X}}X} & {Y\times_{\mathcal{Y}}Y} \\
	X & Y
	\arrow["{f_2}", from=1-1, to=1-2]
	\arrow["{p_2}", shift left=2, from=1-1, to=2-1]
	\arrow["{p_1}"', shift right=2, from=1-1, to=2-1]
	\arrow["{p_2}", shift left=2, from=1-2, to=2-2]
	\arrow["{p_1}"', shift right=2, from=1-2, to=2-2]
	\arrow["{f_1}"', from=2-1, to=2-2]
\end{tikzcd}.\]
     As we assumed \(H\subset \mathcal{I}_f\), for all \((x,y,\alpha)\in X\times_{\mathcal{X}}X(T)\), i.e. \(\alpha:p(x)\cong p(y)\) in \(\mathcal{X}\), and for all \(h\in \text{Aut}_T(p(x))\), 
    \[f_2(x,y,\alpha)=f_2(x,y,\alpha\circ h)\]
    so \(f_2\) descends to a morphism \(\overbar{f}_2:X\times_{\mathcal{X}\!\!\!\fatslash H}X\to Y\times_{\mathcal{Y}}Y\). This induces a factorization of morphisms of groupoids

    \[\begin{tikzcd}
	{X\times_{\mathcal{X}}X} & {X\times_{\mathcal{X}\!\!\!\fatslash H}X} & Y\times_{\mathcal{Y}}Y \\
	X & X & Y
	\arrow["{p_2}", shift left=2, from=1-1, to=2-1]
	\arrow["{p_1}"', shift right=2, from=1-1, to=2-1]
	\arrow[from=1-1, to=1-2]
	\arrow["{\text{id}}"', from=2-1, to=2-2]
	\arrow["{p_2}", shift left=2, from=1-2, to=2-2]
	\arrow["{p_1}"', shift right=2, from=1-2, to=2-2]
	\arrow[shift left=2, from=1-3, to=2-3]
	\arrow[shift right=2, from=1-3, to=2-3]
	\arrow["{\overbar{f}_2}", from=1-2, to=1-3]
	\arrow["{f_1}"', from=2-2, to=2-3]
 \arrow["{p_2}", shift left=2, from=1-3, to=2-3]
	\arrow["{p_1}"', shift right=2, from=1-3, to=2-3]
\end{tikzcd}\]
which descends to the desired factorization of morphisms of stacks, by applying the functoriality of quotient stacks again.

Suppose \(f_1:\mathcal{X}\!\!\!\fatslash H\to \mathcal{Y}\) is another morphism such that \(f=f_1\circ \rho\). Then 
\[f^{rig}\circ \rho \circ p=f\circ p=f_1\circ\rho \circ p\]
and morphisms of groupoids in algebraic spaces induced by \(f^{rig}\) and \(f_1\) are the same, \(f^{rig}=f_1\) proving uniqueness. 
\end{proof}

\begin{corollary}\label{formal corollary}
    Let
\[\begin{tikzcd}
	{\mathcal{X}_1} & {\mathcal{X}_2} \\
	{\mathcal{X}_3} & {\mathcal{X}_4}
	\arrow["{f_{12}}", from=1-1, to=1-2]
	\arrow[""{name=0, anchor=center, inner sep=0}, "{f_{24}}", from=1-2, to=2-2]
	\arrow["{f_{13}}"', from=1-1, to=2-1]
	\arrow[""{name=1, anchor=center, inner sep=0}, "{f_{34}}"', from=2-1, to=2-2]
	\arrow["\alpha"', shorten <=4pt, shorten >=4pt, Rightarrow, from=0, to=1]
\end{tikzcd}\]
    be a \(2\)-commutative diagram of algebraic stacks. Suppose \(H_i\subset \mathcal{I}_{\mathcal{X}_i}\) are flat subgroups of finite presentation over \(\mathcal{X}_i\) and let \(\rho_i:\mathcal{X}_i\to \mathcal{X}_i\!\!\fatslash H_i\) be the associated rigidification morphisms. If \(H_i\subset \mathcal{I}_{\rho_j\circ f_{ij}}\) for all \(1\leq i<j\leq 4\), then the diagram 
    \[\begin{tikzcd}
	{\mathcal{X}_1\!\!\!\fatslash H_1} && {\mathcal{X}_2\!\!\!\fatslash H_2} \\
	{\mathcal{X}_3\!\!\!\fatslash H_3} && {\mathcal{X}_4\!\!\!\fatslash H_4}
	\arrow["{(\rho_2\circ f_{12})^{rig}}", from=1-1, to=1-3]
	\arrow[""{name=0, anchor=center, inner sep=0}, "{(\rho_4\circ f_{24})^{rig}}", from=1-3, to=2-3]
	\arrow["{(\rho_3\circ f_{13})^{rig}}"', from=1-1, to=2-1]
	\arrow[""{name=1, anchor=center, inner sep=0}, "{(\rho_4\circ f_{34})^{rig}}"', from=2-1, to=2-3]
	\arrow["\rho_4(\alpha)"', shorten <=6pt, shorten >=6pt, Rightarrow, from=0, to=1]
\end{tikzcd}\]
is \(2\)-commutative as well.
\end{corollary}
\begin{proof}
    By the universal property, it is enough to show that 
    \[(\rho_4\circ f_{24})^{rig}\circ(\rho_2\circ f_{12})^{rig}\circ \rho_1\cong (\rho_4\circ f_{34})^{rig}\circ(\rho_3\circ f_{13})^{rig}\circ \rho_1 \]
    via \(\rho_4(\alpha)\). Applying the universal property to the left hand side of the equation, one gets
    \[(\rho_4\circ f_{24})^{rig}\circ(\rho_2\circ f_{12})^{rig}\circ \rho_1=(\rho_4\circ f_{24})^{rig}\circ \rho_2\circ f_{12}=\rho_4\circ f_{24}\circ f_{12}\]
    which is isomorphic via \(\rho(\alpha)\) to 
    \[\rho_4\circ f_{34}\circ f_{13}\]
    which is equal to 
    \[ (\rho_4\circ f_{34})^{rig}\circ(\rho_3\circ f_{13})^{rig}\circ \rho_1 \]
    by the universal property.
\end{proof}
\end{appendix}
\begingroup
\sloppy
\RaggedRight
\printbibliography
\endgroup
\end{document}